\UseRawInputEncoding
\documentclass[11pt, a4paper,leqno]{amsart}
\usepackage{amsmath,amsthm,amscd,amssymb,amsfonts, amsbsy}
\usepackage{latexsym}
\usepackage{exscale}
\usepackage{txfonts}
\usepackage[colorlinks,citecolor=red,pagebackref,hypertexnames=false]{hyperref}
\usepackage{enumitem} 
\usepackage{pgf}
\usepackage{color}

\numberwithin{equation}{section}

\theoremstyle{plain}
\newtheorem{theorem}[equation]{Theorem}
\newtheorem{lemma}[equation]{Lemma}

\theoremstyle{definition}
\newtheorem{definition}[equation]{Definition}

\theoremstyle{remark}
\newtheorem{remark}[equation]{Remark}

\newtheorem{claim}[equation]{Claim}

\newcommand{\dv}{\operatorname{div}}

\newcommand{\supp}{\operatorname{supp}}
\newcommand{\dist}{\operatorname{dist}}

\newcommand{\tQ}{Q_*}
\newcommand{\tD}{\Delta_*}

\newcommand{\re}{\mathbb{R}}
\newcommand{\rn}{\mathbb{R}^n}

\newcommand{\ree}{\mathbb{R}^{n+1}}

\newcommand{\eps}{\varepsilon}

\newcommand{\vp}{\varphi}

\newcommand{\p}{\mathcal{P}}

\newcommand{\ppo}{\mathcal{P}\Omega}
\newcommand{\eo}{\partial_e\Omega}
\newcommand{\no}{\partial_n\Omega}
\newcommand{\so}{\partial_s\Omega}
\newcommand{\sso}{\partial_{ss}\Omega}
\newcommand{\ao}{\partial_a\Omega}
\newcommand{\bo}{\mathcal{B}\Omega}
\newcommand{\So}{\mathcal{S}\Omega}

\newcommand{\F}{\mathcal{F}}
\newcommand{\tF}{\mathcal{F}_*}
\newcommand{\Hc}{\mathcal{H}}
\newcommand{\C}{\mathcal{C}({\bf Y_0},\!{\bf Y})}

\newcommand{\m}{\mathcal{M}}
\newcommand{\W}{\mathcal{W}}
\newcommand{\po}{{\partial\Omega}}

\newcommand{\RRR}{\mathcal{R}}

\newcommand{\bX}{{\bf X}}
\newcommand{\bY}{{\bf Y}}
\newcommand{\bx}{{\bf x}}
\newcommand{\by}{{\bf y}}

\newcommand{\hm}{\omega}

\newcommand{\pom}{\partial\Omega}

\DeclareMathOperator{\diam}{diam}

\def\Yint#1{\mathchoice
    {\YYint\displaystyle\textstyle{#1}}%
    {\YYint\textstyle\scriptstyle{#1}}%
    {\YYint\scriptstyle\scriptscriptstyle{#1}}%
    {\YYint\scriptscriptstyle\scriptscriptstyle{#1}}%
      \!\iint}
\def\YYint#1#2#3{{\setbox0=\hbox{$#1{#2#3}{\iint}$}
    \vcenter{\hbox{$#2#3$}}\kern-.51\wd0}}
\def\longdash{{-}\mkern-2.0mu{-}} 
  
\def\tiltlongdash{\rotatebox[origin=c]{12}{$\longdash$}}

\def\tiltfiint{\Yint\tiltlongdash}

\begin{document}

\title[Weak Reverse H{\"o}lder Inequality]{A weak reverse H{\"o}lder inequality for caloric  measure}

\author{Alyssa Genschaw}

\address{Alyssa Genschaw
\\
Department of Mathematics
\\
University of Missouri
\\
Columbia, MO 65211, USA} \email{adcvd3@mail.missouri.edu}

\author{Steve Hofmann}

\address{Steve Hofmann\\
Department of Mathematics\\
University of Missouri\\
Columbia, MO 65211, USA} \email{hofmanns@missouri.edu}

\thanks{The authors were supported by 
NSF grant number DMS-1664047.}

\date{\today}
\subjclass[2000]{42B99, 42B25, 42B37, 35K05, 35K20, 31B20, 31B25}

\keywords{Dirichlet problem, caloric measure, parabolic measure, heat equation,
divergence form parabolic equations, 
weak-$A_\infty$, reverse H\"older inequality, Ahlfors-David 
Regularity 
}

\begin{abstract} 
Following a result of Bennewitz-Lewis for non-doubling harmonic measure,
we prove a criterion for non-doubling caloric measure to satisfy a 
weak reverse H{\"o}lder inequality on an open set $\Omega \subset \ree$,
assuming as a background hypothesis only that the essential boundary
of $\Omega$ satisfies an appropriate parabolic version
 of Ahlfors-David
regularity (which entails some backwards in time thickness).
 We also show that the weak reverse H\"older estimate is equivalent to solvability of
 the initial Dirichlet problem with ``lateral" data in $L^p$, for some $p<\infty$,
 in this setting.
\end{abstract}

\maketitle

{\small
\tableofcontents}

\section{Introduction}\label{s1}  

It is well known that for a Lipschitz domain
$\Omega$, the Dirichlet problem  for a divergence form uniformly
elliptic equation $Lu= -\dv A \nabla u = 0$, with data in $L^p(\partial\Omega)$, is solvable for some $1<p<\infty$ if and only if  elliptic-harmonic measure for $L$ 
is absolutely continuous with respect to surface measure and the Poisson kernel satisfies a reverse H{\"o}lder condition with exponent $p'$; see \cite{Ke} and the references cited there. 
In particular, in the case that $L$ is the Laplacian,  the Poisson kernel satisfies an
$L^2$ reverse H\"older inequality, and therefore the Dirichlet problem is solvable with 
data in $L^2(\pom)$ (see \cite{Da}).

In this paper we prove a parabolic version of a result of 
Bennewitz-Lewis \cite{BL}, who gave a 
criterion for nondoubling harmonic measure to satisfy a weak-$A_\infty$ condition, or
equivalently, for the Poisson kernel to
satisfy a weak reverse H\"older condition; see Definition \ref{defAinfty} below.
To put this work in context, we recall that
David-Jerison \cite{DJ} and Semmes \cite{S} proved that 
harmonic measure $\omega$ on the boundary of an NTA domain with 
Ahlfors-David regular boundary is $A_{\infty}$ with respect to surface measure.   
 The idea of the approach in \cite{DJ}
is  first to prove a geometric result, whereby domains satisfying a  certain
two sided interior and exterior thickness condition (that is, the two sided ``Corkscrew" condition), 
and having ADR boundaries,
could be approximated in a ``Big Pieces" sense by Lipschitz sub-domains.  As a consequence, by the maximum principle
combined with the fundamental result of \cite{Da}, one obtains a certain local ampleness property of the 
harmonic measure (see \eqref{eq1.4} below for the parabolic version), 
which may then, in the presence of the Harnack chain condition, 
be self-improved to give the $A_\infty$ property.

In \cite{BL}, the authors show that this self-improvement procedure, i.e., the passage from local ampleness of harmonic measure to quantitative absolute continuity,
can still be executed, even in the absence of the Harnack chain condition, and as a consequence are able to extend the result of \cite{DJ} and \cite{S}, in an appropriate way,
to much more general domains.  They are able to conclude only  
that harmonic measure is weak-$A_{\infty}$ with respect to surface measure, but on the other hand,
this conclusion is best possible:
their results apply to domains in which harmonic measure need not be doubling
(in particular, to the case that the domain satisfies a uniform interior big pieces of Lipschitz graph condition and an interior corkscrew condition, but no connectivity property, such as the Harnack chain condition).


  The goal of the present paper is to extend the results of \cite{BL} to the parabolic setting.
 As regards geometric hypotheses,
we assume only that $\Omega \subset \ree$ is an open set whose boundary satisfies an appropriate version
of a parabolic Ahlfors-David regularity condition.
In particular, we impose no connectivity hypothesis, such as a parabolic Harnack chain condition.
We may then consider the initial-Dirichlet problem with ``lateral" data in $L^p$, in subdomains of the form
$\Omega^T = \Omega \cap \{t>T\}$, for appropriate fixed times $T$.  We shall return to the latter point below.

We shall consider 
the heat operator
\begin{equation}
\label{eq1.1a}
L_0:=\partial_t- \nabla\cdot\nabla,
\end{equation}
where $ \nabla\cdot\nabla$ is the usual Laplacian in $\rn$, acting in the space variables.
With a caveat, to be discussed momentarily, 
our results may apply more generally to divergence form parabolic operators 
\begin{equation}
\label{eq1.1}
L:=\partial_t- \dv A(X,t) \nabla,
\end{equation}
defined in an open set $\Omega\subset\mathbb{R}^{n+1}$  as described above,
where $A$ is $n\times n$, real,  $L^\infty$, and satisfies the uniform ellipticity condition
\begin{equation}
\label{eq1.1*} \lambda|\xi|^{2}\leq\,\langle A(X,t)\xi,\xi\rangle := \sum_{i,j=1}^{n+1}A_{ij}(X,t)\xi_{j} \xi_{i}, \quad
  \Vert A\Vert_{L^{\infty}(\mathbb{R}^{n})}\leq\lambda^{-1},
\end{equation}
 for some $\lambda>0$, and for all $\xi\in\rn$, and a.e. $(X,t)\in \Omega$.  
We do not require that the matrix $A(X,t)$ be symmetric. 
We reference the paper by Moser \cite{M}, where 
the results are stated under an assumption of
symmetry,  but in fact symmetry is not needed: 
see  \cite{DK}, \cite{SSSZ}, \cite{QX}.

 Some comments are in order.   As mentioned above, there is
 a caveat when applying our results to
 variable coefficient operators, namely that at present it appears to be an open problem to 
 construct parabolic measure for such operators, 
 in the very general class of domains that we consider here.
To do so first requires that one can solve, in an appropriate sense, 
the Dirichlet 
problem with continuous data, so that parabolic measure can be
 constructed via Riesz representation.  
One can construct Perron solutions (as a supremum of subsolutions) 
but then, to apply the Riesz representation theorem, one needs linearity of the solution map
(i.e., the map that sends data $f$ to the value of the solution at a given point $(X,t)$ in the domain).  
For the heat equation this works, since it is known that
continuous functions on the parabolic boundary (or, to be more
precise, on the ``essential boundary"; see Definition \ref{parabolicb}) are resolutive for the 
heat equation (see \cite{W1} or \cite{W2}), and therefore the solution map is linear.  
On the other hand, 
for more general parabolic operators, linearity of the solution map 
would follow if one could solve the continuous Dirichlet problem, in the sense of Definition \ref{bvpdef} below.
A rather general result in this direction was obtained in \cite{CDK}, where the authors
assume an exterior measure condition, backwards in time (see \cite[Definition 1.3]{CDK}).
Otherwise, it would suffice to have a Wiener criterion to ensure continuity up to the parabolic boundary, along with
enough solutions in a class to which the Wiener criterion can be applied; 
to our knowledge, there are versions of the parabolic Wiener test that apply 
to Perron solutions either for the heat equation \cite{La}, \cite{EG}, or to divergence form
parabolic equations with $C^\infty$ or
$C^1$-Dini coefficients \cite{GL}, \cite{FGL}, respectively;  or, in the case of general divergence form parabolic
operators, to some class of weak solutions (either the class $V_2$ \cite{BiM}, or  $W^{1,2}$ \cite{GZ}).  
It appears to be an open problem to construct solutions of the latter sort, say 
for data that is Lipschitz with compact support, except in cylindrical domains \cite{LSU}, in Lip(1,1/2) domains
\cite{BHL}, and in parabolic Reifenberg flat domains \cite{BW}.   
Thus, our results will apply without further qualification to the heat equation, or to operators with
$C^1$-Dini coefficients,
but at present, they
will apply to general divergence form parabolic operators only if one is given {\em a priori} 
that the classical
Dirichlet problem, with continuous data, is solvable.  We observe that the capacitary conditions in \cite{BiM,GZ}
hold in our setting:  they follow from the time-backwards version of ADR (Definition \ref{defbackadr}) 
that we assume; the obstacle to our applying these Wiener criteria, is the lack of solutions. 

Before stating our main theorem, we briefly
introduce some of the concepts and notation
to be used.  All additional terminology used in 
the statement of the theorem, and not  discussed here or above, 
will be defined precisely in the sequel.  For now, we note that all distances and diameters
are taken with respect to the parabolic distance \eqref{pardist}, and that $\delta(X,t) := \dist((X,t),\eo)$,
where $\eo$ denotes the {\em essential boundary} 
(see Definition \ref{parabolicb} below) of an open set $\Omega\subset \ree$.  
We further note that ``surface measure" 
$\sigma$ on the {\em quasi-lateral boundary}\footnote{This comprises all but the initial
part of the essential boundary, and all but the terminal part of the singular boundary; see Definition \ref{parabolicb}.
It may seem more natural to expect that
``surface measure" should be defined on the lateral boundary, rather than on
the quasi-lateral boundary;
in the present work, the ADR condition that we impose will imply that in some sense,
the non-lateral parts of the quasi-lateral boundary are fairly negligible: in particular the singular boundary will exist
only at the terminal time of $\Omega$, hence the quasi-lateral boundary will be a subset of the essential boundary; 
the quasi-lateral 
 boundary then becomes a natural substitute for the lateral boundary.  We 
shall return to this point below:  see Definitions \ref{parabolicb}, \ref{defadr}, and \ref{defbackadr}, and Remarks
\ref{r-adr} and \ref{r-adr2}.} 
$\Sigma$, is defined by
$d\sigma=d\sigma_{\!s}ds$, 
where $d\sigma_{\!s}:=\Hc^{n-1}\vert_{\Sigma_s}$, 
the restriction of 
$(n-1)$-dimensional Hausdorff measure to the time slice $\Sigma_s := \Sigma\cap \{t\equiv s\}$.

We note that for an arbitrary open set $\Omega\subset \ree$, caloric measure
may be constructed via the PWB method, since continuous functions on the essential boundary
are resolutive;  see \cite{W1} or \cite[Chapter 8]{W2}.
 
For a sufficiently large (and eventually fixed)  constant $K_1$, given $(X,t) \in \Omega$, set
\begin{equation}\label{eq1.3}
Q_{X,t} := Q\big((X,t), K_1 \delta(X,t)\big)\,, \quad \Delta_{X,t}:= Q_{X,t}\cap\Sigma\,,
\end{equation}
where in general the parabolic cube $Q\big((X,t), r\big)$ is defined as in \eqref{eqcubedef} below.
Eventually, we shall fix $K_1$ in \eqref{eq1.3} 
large enough, depending
only on the constants in Lemma \ref{Bourgain}.  

For $(X,t) \in \Omega$, and $\Delta_{X,t}$  defined as in \eqref{eq1.3}, 
we shall say that caloric (or parabolic) measure $\omega^{X,t}$ is {\em locally ample} on $\Delta_{X,t}$,
or more precisely,
{\em $(\theta,\beta)$-locally ample},  if
there exists  constants $\theta,\,\beta\in (0,1)$ such that 
\begin{align}\label{eq1.4}
\sigma(E)\geq 
(1-\theta)\sigma(\Delta_{X,t})\,\,\implies\,\,
\omega^{X,t}(E)= \omega_L^{X,t}(E)\geq\beta\,,
\end{align}
where $E\subset \Delta_{X,t}$ is a Borel set.  We observe that if \eqref{eq1.4} holds for some $K_1=K\geq 2$, 
then it also
holds with $K_1 =K'>K$, for some $\theta'=\theta'(\theta, K,K',n,ADR)$.  Thus, we may always fix a larger value of
$K_1$, at our convenience.

Set 
$ T_{min}:=\inf\{T: \Omega\cap\{t\equiv T\} \neq \emptyset\}$ (note that 
we may have $T_{min} =-\infty$).

The main result of this paper is the following.  
 Precise definitions of terminology may be found in the sequel.


\begin{theorem}\label{tmain}   Let $\Omega\subset \ree$ be 
an open set 
whose quasi-lateral boundary $\Sigma$ is globally ADR. 
Let $(x_0,t_0)\in\Sigma$, and let 
$0<r< 
\sqrt{t_0-T_{min}}/(8\!\sqrt{n})$. 
Assume that 
$\Sigma$ is time-backwards ADR 
on $ \Delta_{2r} =\Sigma \cap Q_{2r}(x_0,t_0)$, 
and suppose that there are constants $\theta,\,\beta\in (0,1)$, and a value of $K_1\geq 2$ in \eqref{eq1.3},
 such that
caloric measure $\omega^{X,t}$ satisfies the $(\theta,\beta)$-local ampleness condition 
\eqref{eq1.4} on $\Delta_{X,t}$ for each
$(X,t) \in \Omega \cap Q_{2r}(x_0,t_0)$. 

Then there exist constants $C\geq 1$, $\gamma >0$,  
such that if $(Y_0,s_0) \in \Omega \setminus Q_{4r}(x_0,t_0),$ 
then $\omega^{Y_0,s_0} \ll\sigma$ on $\Sigma \cap Q_r(x_0,t_0)$, with
$d\omega^{Y_0,s_0}/d\sigma=h$ satisfying
\begin{multline}\label{weakRH}
\left(\rho^{-n-1}\iint_{\Delta_{\rho}(y,s)}h^{1+\gamma}d\sigma\right)^{1/(1+\gamma)} \leq 
C\rho^{-n-1}
\iint_{\Delta_{2\rho}(y,s)}h\, d\sigma 
\\=C\rho^{-n-1}\omega^{Y_0,s_0}\left(\Delta_{2\rho}(y,s)\right),
\end{multline}
whenever $(y,s)\in\Sigma$ and $Q_{2\rho}(y,s)\subset Q_{r}(x_0,t_0)$, where 
$\Delta_{\rho}(y,s)=Q_{\rho}(y,s)\cap \Sigma$, and 
$\Delta_{2\rho}(y,s)=Q_{2\rho}(y,s)\cap \Sigma$.

Moreover, the same result applies to the parabolic measure associated to
a uniformly parabolic divergence form operator $L$, provided 
that the continuous 
Dirichlet  problem is solvable for $L$ in $\Omega$ (see Definition \ref{bvpdef})
(in particular, this is true if the coefficients are $C^1$-Dini).
\end{theorem}
To clarify matters, we remark that by the ADR hypothesis on $\Sigma$ (see Definition \ref{defadr}),
we have that $\rho^{n+1} \approx \sigma(\Delta_{\rho}(y,s)) \approx \sigma(\Delta_{2\rho}(y,s))$.
The time-backwards ADR condition is an enhanced version of ADR, which entails some thickness of
$\Sigma$ backwards in time; see Definition \ref{defbackadr}, and Remarks \ref{r1.22}, \ref{r-adr2},
and \ref{1.8}.

Similar results in the parabolic setting have previously been established 
under the more restrictive assumptions that 1) $L$ is the heat operator and the lateral boundary
of the domain is given locally as the graph 
of a function $\psi(x,t)$ which is Lipschitz in the space variable, and has a 1/2-order time derivative in parabolic 
BMO  \cite{LM}\footnote{For domains whose
 lateral boundary is given locally as a graph, the 1/2 order derivative in BMO condition of 
 \cite{LM} is in the nature of best possible:
 there is a counterexample of Kaufmann and Wu \cite{KW}, with $\psi \in$ 
Lip$_{1/2}$ in the time variable.}, 2)
$\Omega = \{(x_0,x,t) \in (0,\infty)\times \re^{n-1}\times \re\}$ is a half-space and the coefficients
of $L$ satisfy a certain Carleson measure regularity property \cite{HL},  and 3) $L$ is the heat operator 
and either $\Omega$ is a parabolic 
Reifenberg flat domain \cite{HLN}, or  $\Omega$ is a parabolic chord-arc domain  \cite{NS}; in each of these
settings, $\Omega$ enjoys a parabolic version of the Harnack Chain condition, which entails a rather strong
quantitative version of connectivity.  As mentioned above,
the elliptic analogue of our result was proved in \cite{BL}, without any connectivity hypothesis.
The new contribution of the present paper is to dispense with all connectivity assumptions,
both qualitative and quantitative, in the parabolic setting.  
The elliptic version obtained in \cite{BL} has proved to be useful
in various applications, see, e.g., \cite{HLe} and \cite{HM}.  We shall discuss 
two applications of our work in the sequel (see Section \ref{s5a}).

The paper is organized as follows.  In the remainder of this section, 
we present some basic notations and definitions.  In Section \ref{s2}, we state four lemmas which we then use to 
prove Theorem \ref{tmain}; we also
state Theorem \ref{t2.9}, concerning the equivalence between the weak-$A_\infty$ property and $L^p$ solvability of the 
initial-Dirichlet problem. In Section \ref{s3} we prove Lemma \ref{lemma1}, and in Section \ref{q} we 
prove Theorem \ref{t2.9}. In Appendix \ref{s4} we prove a Bourgain-type estimate (Lemma \ref{Bourgain}), 
and in Appendix \ref{s5} we prove H{\"o}lder continuity at the boundary (Lemma \ref{continuity}).
In Appendix \ref{s6}, we prove a technical fact about the essential boundary (Lemma \ref{essbclosed}).

\medskip

\noindent{\bf Notation and Definitions}. 
For a set $A\subset \ree$, we define
\begin{equation}\label{tminmax}
T_{min}(A):=\inf\{T: A\cap\{t\equiv T\} \neq \emptyset\}\,,\quad T_{max}(A):=\sup\{T: A\cap\{t\equiv T\} \neq \emptyset\}
\end{equation}
(note:  it may be that $T_{min}(A) = -\infty$, and/or that $T_{max}(A) = +\infty$).  In the special case that $A=\Omega$, an open set that
has been fixed, we will simply write
$T_{min} = T_{min}(\Omega)$ and $T_{max}=T_{max}(\Omega)$.

\begin{definition}[{\bf Parabolic cubes}]\label{paraboliccube}
An (open) parabolic cube in $\rn \times \re$ with  center
$(X,t)$:
\begin{multline}\label{eqcubedef}
Q_r(X,t):=Q((X,t),r)\\:=\{(Y,s)\in \rn \times \re: |X_i-Y_i|<r \,,  \, 1\leq i \leq n, \,
t-r^2<s < t+r^2\}.
\end{multline}
With a mild abuse of terminology, we refer to $r$ as the ``parabolic sidelength" (or simply the ``length")
of $Q_r(X,t)$.
We shall sometimes simply write $Q_r$ to denote a cube of parabolic length $r$, when the center is implicit,
and for $Q=Q_r$, we shall write $\ell(Q) = r$.

We also consider the time-backward and time-forward versions:  
\begin{multline*}
Q^-((X,t),r):=Q_r^-(X,t)\\:=\{(Y,s)\in \rn\times \re: |X_i-Y_i|<r\,,\, 
1\leq i\leq n\,, \, t-r^2<s< t\},
\end{multline*}
\begin{multline*}
Q^+((X,t),r):=Q_r^+(X,t)\\:=\{(Y,s)\in \rn\times \re : |X_i-Y_i|<r \,,\,
1\leq i\leq n\,,\, t< s < t+r^2\}\,.
\end{multline*}
\end{definition}

\begin{definition}[{\bf Classification of boundary points}]\label{parabolicb}
Following \cite{L}, given an open set $\Omega\subset \ree$, we define its {\em parabolic
boundary} $\p\Omega$ as
$$\mathcal{P}\Omega:=\left\{(x,t)\in\partial\Omega: \forall r>0 \,, \,Q_r^-(x,t)\,
\text{ meets }\, \ree\setminus \Omega\right\}.
$$
 The {\em bottom boundary}, denoted $\bo$, is defined as 
 \[\bo:= \left\{(x,t)\in\ppo: \exists \,\eps>0 \,  \text{ such that } \,Q_\eps^+(x,t)\,
\subset \Omega\right\}.\]
The {\em lateral boundary}, denoted $\So$, is defined as
$\So :=\ppo\setminus \bo$.

 Following \cite{W1,W2}, we also define the {\em normal boundary}, denoted
 $\no$, to be equal to the parabolic boundary in a bounded domain, while in an unbounded domain, we 
 append the point at infinity:
 $\no = \ppo\cup\{\infty\}$.   The {\em abnormal boundary} is defined as $\ao:= \pom\setminus \no$, thus:
  \[\ao:= \left\{(x,t)\in\partial\Omega: \exists \,\eps>0 \,  \text{ such that } \,Q_\eps^-(x,t)\,
\subset \Omega\right\}.\]
The abnormal boundary is further decomposed into $\ao=\so\cup\sso$ (the
{\em singular boundary} and {\em semi-singular boundary}, respectively), where
\[\so:=  \left\{(x,t)\in\partial_a\Omega: \exists \,\eps>0 \,  \text{ such that } \,Q_\eps^+(x,t)\,
\cap \Omega=\emptyset\right\},\]
and
\[\sso:=  \left\{(x,t)\in\partial_a\Omega: \forall \,r>0 \,  \,Q_r^+(x,t)\,\text{ meets }\,
\Omega\right\}.\]
The {\em essential boundary} 
$\eo$, is defined as
\begin{equation}\label{sigma}
\eo\,:=\, \no\,\cup\, \sso \,=\, \partial\Omega \setminus \so
\end{equation}
(where we replace $\pom$ by $\pom\cup\{\infty\}$ if $\Omega$ is unbounded).
Finally, 
 we define the {\em quasi-lateral boundary} $\Sigma$ to be
\begin{equation}\label{sigmazero}
\Sigma:= \,\left\{
\begin{array}{l}
\po \,,\quad  \text{if } T_{min} = -\infty \, \text{ and } T_{max} = \infty\\[4pt]
\po \setminus (\bo)_{T_{min}}\,, \quad  \text{if } T_{min} > -\infty \, \text{ and } T_{max} = \infty\\[4pt]
\po\setminus (\so)_{T_{max}} \,, \quad  \text{if } T_{max} <\infty \, \text{ and } T_{min} = - \infty
 \\[4pt]
 \po\setminus \left( (\bo)_{T_{min}} \cup (\so)_{T_{max}}\right) \,, 
 \quad  \text{if }\, -\infty<T_{min}<T_{max} <\infty \,.
 \end{array}\right.
\end{equation}
where $(\bo)_{T_{min}}$ is the 
time slice of $\bo$ with $t\equiv T_{min}$,
and $(\so)_{T_{max}}$ is the time slice of $\so$ with $t \equiv T_{max}$. 
Observe that for a cylindrical domain
$\Omega= U \times (T_{min},T_{max})$, with $U\subset \rn$ a domain in the spatial variables,  
then $\Sigma$ would simply be the usual lateral boundary.
\end{definition}
 Caloric measure is supported on the essential boundary; see
 \cite{Su}, or  \cite{W1,W2}.


For future reference, we record here the following fact.
\begin{lemma}\label{essbclosed} The essential boundary $\eo$, and the quasi-lateral boundary 
$\Sigma$, are closed sets.
\end{lemma}

We defer the proof of this lemma to Appendix \ref{s6}.




\begin{itemize}

\item We use the letters $c,C$ to denote harmless positive constants, not necessarily
the same at each occurrence, which depend only on dimension and the
constants appearing in the hypotheses of the theorems (which we refer to as the
``allowable parameters'').  We shall also
sometimes write $a\lesssim b$ and $a \approx b$ to mean, respectively,
that $a \leq C b$ and $0< c \leq a/b\leq C$, where the constants $c$ and $C$ are as above, unless
explicitly noted to the contrary.  
\item 
We shall
use lower case letters $x,y,z$, etc., to denote the spatial component of points on the boundary $\pom$, 
and capital letters
$X,Y,Z$, etc., to denote the spatial component of generic points in $\ree$ (in particular those in $\Omega$).

\item For the sake of notational brevity, we shall sometimes also use boldface capital letters
to denote points in space time $\ree$, and lower case boldface letters to denote points on $\pom$;  thus,
\[ \bX = (X,t) , \quad \bY = (Y,s) , \quad
\text{and} \quad  \bx = (x,t) , \quad \by = (y,s) ,\]


\item We shall orient our coordinate axes so that time runs from left to right.
\item $\Hc^{d}$ denotes $d$-dimensional Hausdorff measure.
\item For $A\subset \ree$, let $A_{s}:=\{(X,t)\in A: t\equiv s\}$ denote the time slice of $A$ with $t\equiv s$. 
\item We  let
$d\sigma=d\sigma_{\!s}ds$ denote the ``surface measure'' on the quasi-lateral boundary
$\Sigma$, where $d\sigma_{\!s}:=\Hc^{n-1}\vert_{\Sigma_s}$, and $\Sigma_s$ is the time slice of $\Sigma$, with
$t\equiv s$.  See Remark \ref{r-adr} for some clarifying comments.
\item The parabolic norm of a vector ${\bf X}\in\ree$ is defined as 
\begin{align}\label{pardist}
\|\bX \|=||(X,t)||=|X|+|t|^{1/2},
\end{align}
and we refer to the distance induced by this norm as the parabolic distance.
\item If $\bX\in \Omega$, we set $\delta(\bX):=\text{dist}(\bX, \eo),$ the 
parabolic distance to the essential boundary.
\item For a set $A\subset \ree$, we shall write $\diam(A)$ to denote the 
diameter of $A$ with respect to the parabolic distance, i.e.,
\begin{equation}\label{pardiam}
\diam(A) := \sup_{\left({\bf X},{\bf Y}\right)\in A\times A}\|{\bf X} - {\bf Y}\|\,.
\end{equation}

\item Given a Borel measure $\mu$, and a Borel set $A\subset \rn$, with positive and finite $\mu$ measure, we
set $\fint_A f d\mu := \mu(A)^{-1} \int_A f d\mu$;  if $A$ is a subset of space-time $\ree$, we then write
$\tiltfiint_A f d\mu := \mu(A)^{-1} \iint_A f(X,t)\, d\mu(X,t)$.

\item A ``surface cube" on $\Sigma$ is defined by
$$\Delta =Q\cap \Sigma\,,$$
where $Q$ is a parabolic cube centered on $\Sigma$,
or more precisely,
$$\Delta = \Delta_r(x,t):= Q_r(x,t)\cap\Sigma\,,$$
with $(x,t)\in \Sigma$.  We note that the
``surface cubes'' are not the same as the dyadic cubes
of M. Christ \cite{Ch}
 on $\Sigma$; we apologize to the reader for the possibly confusing terminology.
 
 \end{itemize}
 
 

\begin{definition}\label{bvpdef}
We define the following boundary value problems.  
The second
is relevant only in the case that $T_{min} = -\infty$.

 \begin{list}{}{\leftmargin=0.4cm  \itemsep=0.2cm}


   \item I. Continuous Dirichlet Problem:   
 \begin{center}
$(D)\left\{ \begin{array}{rl}Lu&\!\!=0  \textrm{ in } \Omega  \\
u\vert_{\eo}&\!\!=f \in C_c(\eo)
\\ u&\!\!\in C(\Omega\cup\no)\,.
\end{array} \right.$
\end{center} 
If $\Omega$ is unbounded, 
we further specify that $\lim_{\|\bX\|\to\infty}
u(\bX)=0$.
Here,  we interpret the statement $u\vert_{\eo}=f \in C_c(\eo)$
to mean that 
 \[\lim_{(X,t)\to (y,s)} u(X,t) = f(y,s)\,,\quad (y,s) \in \no\,,\]
and
 \[\lim_{(X,t)\to (y,s^+)} u(X,t) = f(y,s)\,,\quad (y,s) \in \sso\,.\]
If  the preceeding problem is solvable for all $f\in C_c(\eo)$, then we say that the
``continuous Dirichlet problem is solvable for $L$."  

 
 \item II. $L^p$ Dirichlet Problem: 
 \begin{center}
$(D)_p\left\{ \begin{array}{rl}Lu&=0  \textrm{ in } \Omega  \\
u\vert_{\Sigma}&=f\in L^p(\Sigma)
\\&N_* u\in L^p(\Sigma)\,.
\end{array} \right.$
\end{center} 
 
  \item III.  Continuous Initial-Dirichlet Problem: 
 \begin{center}
$(I\text{-}D)\left\{ \begin{array}{rl}Lu&=0  \textrm{ in } \Omega^T:=\Omega\cap \{t>T\}  \\
u(X,T)&=0  \text{ in } \Omega_{T}=  \Omega\cap \{t\equiv T\} \\
u\vert_{\Sigma^T}&=f\in C_c(\Sigma^T)
\\& u\in C(\Omega^T \cup \no^T)\,.
\end{array} \right.$
\end{center} 
Here, $\Sigma^T$ denotes the quasi-lateral boundary of the domain
$\Omega^T$.
The statement $u\vert_{\Sigma^T}=f\in C_c(\Sigma^T)$ is intepreted as in problem I,
and if  $\Omega^T$ is unbounded, 
we further specify that $\lim_{\|\bX\|\to\infty}
u(\bX)=0$.  
 
 \item IV. $L^p$ Initial-Dirichlet Problem: 
 \begin{center}
$(I\text{-}D)_p\left\{ \begin{array}{rl}Lu&=0  \textrm{ in } \Omega^T:=\Omega\cap \{t>T\}  \\
u(X,T)&=0  \text{ in }  \Omega_{T}=  \Omega\cap \{t\equiv T\}\\
u\vert_{\Sigma^T}&=f\in L^p(\Sigma^T)
\\&N_* u\in L^p(\Sigma^T)\,.
\end{array} \right.$
\end{center} 
\end{list}
In problems II and IV, the statement $u\vert_{\Sigma}=f\in L^p(\Sigma)$ (resp., 
$u\vert_{\Sigma^T}=f\in L^p(\Sigma^T)$) is understood in the sense of
parabolic non-tangential convergence.  We shall discuss this issue, as well as the precise definition of the non-tangential
maximal function $N_*u$, in the sequel. In problems III and IV, the statement
$u(X,T)=0  \text{ in } \Omega_{T}$ means that $u$ vanishes continuously on $\Omega_T$.
\end{definition}

\begin{definition}\label{parabolic}({\bf Caloric and Parabolic Measure})
Let $\Omega\subset \ree$ be an open set. 
Let $u$ be the  
PWB solution (see \cite{W1}, \cite[Chapter 8]{W2}) 
of the Dirichlet problem  for the heat equation, with data $f \in C_c(\eo)$. 
By the Perron construction, 
for each point $(X,t) \in \Omega$, 
the mapping 
$f \mapsto u(X,t)$ is  bounded, and by
the resolutivity of functions $f\in C(\eo)$ (see \cite[Theorem 8.26]{W2}), 
it is also linear.
The caloric 
measure with pole $(X,t)$ is the probability measure $\omega^{X,t}$ given by 
the Riesz representation 
theorem, such that
\begin{equation}\label{parmeasuredef}
u(X,t)=\iint_{\eo}f(y,s)\,d\omega^{X,t}(y,s).
\end{equation}
For a general divergence form parabolic operator  $L$ as in \eqref{eq1.1}-\eqref{eq1.1*}, 
parabolic measure $\hm^{X,t}=\hm_L^{X,t}$ may be defined similarly, provided that  the continuous Dirichlet problem is solvable for $L$.
\end{definition}


\begin{definition}\label{defadr} ({\bf  ADR})  (aka {\it Ahlfors-David regular [in the parabolic sense]}).
Let $\Omega\subset \ree$. We say that  the quasi-lateral  boundary 
$\Sigma$ 
 is {\em globally ADR} (or just ADR)
if  there is a constant $M_0$ such that for every parabolic cube $Q_r = Q_r(x,t)$, 
centered on $\Sigma$, 
and corresponding surface cube
$\Delta_r=Q_r\cap \Sigma$, with $r<\diam(\Omega)$,
\begin{equation}\label{eq.adr}
\frac1{M_0} r^{n+1}\leq \sigma(\Delta_r) \leq M_0 r^{n+1}\,.
\end{equation}
 We also say that $\Sigma$ is ADR on a surface cube $\Delta=Q\cap \Sigma$, if there is
a constant $M_0$ such that \eqref{eq.adr} holds for every surface cube $\Delta_r=Q_r\cap\Sigma$,
with $Q_r \subset Q$ and centered on $\Sigma$.
\end{definition}


\begin{definition}\label{defbackadr} 
({\bf Time-Backwards ADR, aka TBADR})
Given a parabolic cube $Q$ centered on $\Sigma$, 
and corresponding surface cube
$\Delta=Q\cap \Sigma$,
we say that $\Sigma$ is time-backwards ADR on
$\Delta$ if it is ADR on $\Delta$, and if, in addition there exists a constant $b>0$ 
such that
\begin{equation}\label{eq.backadr} 
br^{n+1}\leq \sigma(\Delta^-_r)\,, 
\end{equation}
for every
$\Delta^-_r=Q^-_r\cap \Sigma,$ where
$Q_r\subset Q$ is centered at some 
point $(x,t)\in \Sigma$. 

 If  $\Sigma$ is time-backwards ADR on every
$\Delta =\Sigma \cap Q_r(x_0,t_0)$, for all $(x_0,t_0)\in\Sigma$, and for all $r$ with 
$0<r<\sqrt{t_0-T_{min}}/(4\!\sqrt{n})$,  
then we shall simply say that $\Sigma$ is (globally) time-backwards ADR (and
we shall refer to such $\Delta$ as ``admissible"; note that if $T_{min} = -\infty$, then there is no restriction on $r$, and in that case every surface cube is admissible). 
\end{definition}

\begin{remark} \label{r1.22} The assumption of some backwards in time thickness, 
as in Definition \ref{defbackadr},
is rather typical in the parabolic setting.
 See, e.g., the backwards in time
 capacitary conditions in \cite{La}, \cite{EG}, \cite{GL}, \cite{FGL}, \cite{GZ}, \cite{BiM}. 
Moreover, it is not hard to verify that by the result of \cite{EG} (or of \cite{GL}, \cite{FGL}),  
time-backwards ADR on some surface cube $\Delta$ implies 
parabolic Wiener-type regularity of each point in $\Delta$ (and thus global time-backwards
ADR implies regularity of the parabolic boundary $\ppo$).  
\end{remark}

 \begin{remark}\label{r-adr}
By \cite[Theorem 8.40]{W2}, the abnormal boundary $\ao$ is contained in a 
countable union of hyperplanes orthogonal
to the $t$-axis.  
Moreover, the same is true for the bottom boundary $\bo$, since its image under the change of variable
$t \to -t$ is contained in $\ao^*$, for the domain $\Omega^*$ obtained from $\Omega$ by the 
same change of variable.  Thus, $\sigma(\bo) = 0 =\sigma(\ao)$.
 \end{remark}

\begin{remark}\label{r-adr2} The time-backwards ADR condition rules out pathologies like a 
vertical face (with time running from left to right horizontally) on $\Sigma$. 
In particular,  $\sso =\emptyset=\so\setminus \{t\equiv T_{max}\}$, at least locally on any surface cube $\Delta$
on which TBADR holds, and thus $\eo =\no =\ppo=\Sigma$ on such $\Delta$.  Moreover, under the hypotheses of Theorem \ref{tmain},
$\bo \cap \Delta_r(x_0,t_0)$ is fairly negligible: by Remark \ref{r-adr}, this set has $\sigma$ measure zero, 
and thus by the conclusion of Theorem \ref{tmain}, it also has caloric/parabolic measure zero.
\end{remark}

\begin{remark} \label{r1.eqdist}  The significance of the admissibility constraint is as follows.
Recall that $\delta({\bf Y}) := \dist({\bf Y},\eo)$, where the distance is of course the parabolic distance.
We note that, by elementary geometry and 
 \eqref{sigmazero}  (i.e.,  the definition of $\Sigma$), for $(x_0,t_0)\in \Sigma$, and
 for all $ {\bf Y} \in \Omega \cap Q_{r}(x_0,t_0)$, assuming the global time-backwards ADR property
 and using the observations in Remark \ref{r-adr2},
 we have that
\begin{equation*}
r< \sqrt{t_0-T_{min}}/(4\!\sqrt{n}) \implies
\delta({\bf Y}) = \dist({\bf Y}, \Sigma)\,.
\end{equation*}
\end{remark}

\begin{remark}\label{1.8} We will show, in Claim 1 of Appendix
\ref{s4}, that time-backwards ADR yields an apparently stronger property: specifically,  
we show that if $\Sigma$ 
is time-backwards ADR on $\Delta=\Delta_r=\Sigma\cap Q_r(x_0,t_0)$, 
then \eqref{eq.backadr} continues to hold (with a slightly different constant) 
with $\Delta_r^-$ replaced by 
$\Delta_r^-\cap \{t<t_0-(ar)^2\}$, and hence also by
$\Sigma\cap Q^-_r(x_0,t_0-(ar)^2)$, for some uniform $a\in (0,1/2)$.
\end{remark}

\begin{definition}\label{defAinfty}
({\bf $A_\infty$}, weak-$A_\infty$, and weak-$RH_q$). 
Given a closed parabolic ADR set $E\subset\ree$, 
and a surface cube
$\Delta_0:= Q_0 \cap E$,
we say that a Borel measure $\mu$ defined on $E$ belongs to
weak-$A_\infty(\Delta_0)$ if 
for each surface cube $\Delta = Q\cap E$, with $2Q\subseteq Q_0$,
\begin{equation}\label{eq1.wainfty}
\mu (F) \leq C \left(\frac{\sigma(F)}{\sigma(\Delta)}\right)^\theta\,\mu (2\Delta)\,,
\qquad \mbox{for every Borel set } F\subset \Delta\,.
\end{equation}
We recall that, as is well known, the condition $\mu \in$ weak-$A_\infty(\Delta_0)$
is equivalent to the property that $\mu \ll \sigma$ in $\Delta_0$, and that for some $q>1$, the
Radon-Nikodym derivative $k:= d\mu/d\sigma$ satisfies
the weak reverse H\"older estimate 
\begin{equation}\label{eq1.wRH}
\left(\,\, \tiltfiint_\Delta k^q d\sigma \right)^{1/q} \,\leq C\, \tiltfiint_{2\Delta} k \,d\sigma\, 
\approx\,  \frac{\mu(2\Delta)}{\sigma(\Delta)}\,,
\quad \forall\, \Delta = Q\cap E,\,\, {\rm with} \,\, 2Q\subseteq Q_0\,.
\end{equation}
We shall refer to the inequality in \eqref{eq1.wRH} as
an  ``$RH_q$" estimate, and we shall say that $k\in RH_q(\Delta_0)$ if $k$ satisfies \eqref{eq1.wRH}.

If \eqref{eq1.wainfty} holds with $\mu(\Delta)$ in place of $\mu(2\Delta)$, for all $Q\subset Q_0$,
then we say that $\mu \in A_\infty(\Delta_0)$. 
\end{definition}

\section{Lemmas and Proof of Theorem 1.3}\label{s2}
In this section, we state four lemmas which when combined allow us to prove Theorem 1.3. 
We also state Theorem \ref{t2.9}.    We recall that $\Sigma$ 
denotes the quasi-lateral 
 boundary; 
see Definition \ref{parabolicb} and \eqref{sigmazero}.  In the sequel, $L$ is either the heat operator,
or else a divergence form 
parabolic operator 
for which the continuous Dirichlet problem is solvable in $\Omega$
(in particular, this is true if $L$ has $C^1$-Dini coefficients),  and
$\hm =\hm_L$ is the associated caloric/parabolic measure.

Let $a>0$ be the constant mentioned in Remark \ref{1.8}.  In the sequel, $\Omega$ will always denote
an open set in $\ree$, with quasi-lateral  boundary $\Sigma$.  Given a fixed time
$T<\infty$, we set 
\begin{equation}\label{ETdef}
E(T):= \left\{(X,t)\in\ree:\, t<T\right\}\,.
\end{equation}

\begin{lemma}[Parabolic Bourgain-type Estimate]\label{Bourgain} 
Let $(x_0,t_0)\in\Sigma$, and let $0<r< \sqrt{t_0-T_{min}}/(4\!\sqrt{n})$. 
Assume that 
$\Sigma$ is time-backwards ADR on $ \Delta_r := Q_r(x_0,t_0) \cap \Sigma$. 
Then there exists $M_1,\eta>0$ such that for all $(X,t)\in Q_{\frac{a}{M_1}r}\cap \Omega$, 
\begin{equation*}
\omega^{X,t}(\Delta_r)=\omega^{X,t}\big(\Delta_{r}\cap E(T)\big)\geq \eta\,,
\end{equation*}
where $Q_{\frac{a}{M_1}r}:=Q((x_0,t_0),\frac{a}{M_1}r)$, and 
$T:=T_{max}(Q_{\frac{a}{M_1}r})=t_0 +(aM_1^{-1}r)^2$.
\end{lemma}
\begin{remark}
The proof of the Bourgain-type estimate can be found in Appendix \ref{s4}. We 
remark that this estimate could probably
also be derived using capacitary methods found in \cite{EG}. 
We give a direct proof adapting Bourgain's argument to the parabolic setting.
\end{remark}


\begin{remark}
We also obtain a 
Bourgain-type estimate for supersolutions; 
see Appendix \ref{s5}.
\end{remark}

Given a fixed time $T$, and a cube $Q_r$ centered on $\Sigma$, we 
set $\Omega_r:= Q_r\cap\Omega$, and $\Omega_r(T):= \Omega_r\cap E(T)$, with $E(T)$ 
defined as in \eqref{ETdef}.

The following lemma is a consequence of the supersolution version of Lemma \ref{Bourgain}.
The proof will be given in Appendix \ref{s5}.

\begin{lemma}\label{continuity}
Let $(x_0,t_0)\in \Sigma$, and fix $r$ with
$0<r< \sqrt{t_0-T_{min}}/(8\!\sqrt{n})$.  
Set $Q_r:=  Q_{r}(x_0,t_0)$, $Q_{2r} :=  Q_{2r}(x_0,t_0)$, and
suppose that 
$\Sigma$ is time-backwards ADR on 
$ \Delta_{2r} := 
 Q_{2r} \cap \Sigma$.  
Let $u$ be the
parabolic 
measure solution corresponding to non-negative data $f\in C_c(\eo)$, with 
$f\equiv 0$ on $\Delta_{2r}$. 
Then for some $\alpha>0$,
\begin{align*}
u(Y,t)\leq C\left(\dfrac{\delta(Y,t)}{r}\right)^{\alpha}\dfrac{1}{| Q_{2r}\cap E(T_1)|}
\iint_{ \Omega_{2r}(T_1)} u, \hspace{.1in} \forall (Y,t)\in  \Omega_r, 
\end{align*}
where $T_1:= T_{max}(Q_r) =t_0+r^2$, and where the constants $C$ and $\alpha$ depend only on $n,$ $\lambda$, and 
the ADR and time-backwards ADR constants.
\end{lemma}

We observe that in the sequel, it will suffice to have a slightly 
less sharp version of Lemma \ref{continuity}, in which $\Omega_{2r}(T_1)$ is replaced by the larger set
$\Omega_{2r}$; see the proof of Theorem \ref{t2.9} in Section \ref{q} below.

\begin{remark}\label{r2.5}
We now fix $K_1$ in \eqref{eq1.3} to be $K_1:= 20a^{-1} M_1$, where
$M_1, a$ are the constants from Lemma \ref{Bourgain} (so that in turn, $a$ is the constant  
in Remark \ref{1.8}).  With this choice of $K_1$, we then define
$Q_{X,t}$ and $\Delta_{X,t}$ as in \eqref{eq1.3}.  Our assumption in Theorem \ref{tmain} is that
$(\theta,\beta)$-local ampleness
holds for some value of $K_1\geq 2$.  If this $K_1$ exceeds the value defined above, 
then we simply take $M_1$ larger.  On the other hand, as noted previously, if 
$(\theta,\beta)$-local ampleness holds for a smaller value of $K_1 \geq 2$ than that defined above, then it also holds 
for larger $K_1$ (for a possibly different value of $\theta$).  In any case, we are at liberty to fix
$K_1$ as above.
\end{remark}

\begin{lemma}\label{lemma1} 
 Let $(x_0,t_0)\in \Sigma$, and let $0<r< \sqrt{t_0-T_{min}}/(8\!\sqrt{n})$.
Suppose that $\Sigma$ is time-backwards ADR on $\Delta_{2r} := Q_{2r}(x_0,t_0) \cap \Sigma$.
Suppose further
that there exist constants $\theta$, $\beta\in (0,1)$, such that 
$\omega^{X,t}$ satisfies the $(\theta,\beta)$-local ampleness condition 
\eqref{eq1.4} on $\Delta_{X,t}$ for each
$(X,t) \in \Omega \cap Q_{2r}(x_0,t_0)$. 


Then given $\epsilon >0$, there exists $\eta=\eta(\epsilon, n)$, $0<\eta<1$ and $C_{\epsilon}=C(\epsilon,n)\geq 1$ such that for any Borel set
$F\subset \Delta_{2r}$, 
with $\sigma(F)\geq (1-\eta)\sigma ( \Delta_{2r})$, we have
\begin{align}\label{eq2.4}
\omega^{Y,s}(\Delta_r(x_0,t_0))\leq 
\epsilon \omega^{Y,s}(\Delta_{2r}(x_0,t_0))+ C_{\epsilon}\omega^{Y,s}(F),
\end{align}
whenever $(Y,s)\in \Omega\setminus Q_{4r}(x_0,t_0)$. 
\end{lemma}

Lemma \ref{lemma1} is in some sense the main result of this paper.  Along with the next lemma,
it underlies the proof of Theorem \ref{tmain}.

\begin{lemma}\label{lemma2} Let 
$\Sigma$ be a closed ADR set with constant $M_0$. Let $\mu$ be a positive Borel measure on $\ree$ with support contained in $\Sigma$, and $\mu(\Sigma)\leq1$. 
Suppose for some $(x,t)\in \Sigma$, $r>0$, there exists 
$\epsilon_1, \zeta>0, C_1\geq 1$ such that
\begin{align*}
\hspace{.1in}\mu(Q_{\rho}(z,\tau))
\leq \epsilon_1 \mu(Q_{2\rho}(z,\tau))+C_1\mu(P),
\end{align*}
whenever $P\subset \Delta_{2\rho}(z,\tau):= \Sigma
\cap Q_{2\rho}(z,\tau)$ is a Borel set with
\begin{align*}
\sigma(P) \geq (1-\zeta) \,\sigma( \Delta_{2\rho}(z,\tau))
\end{align*}
and $(z,\tau)\in\Sigma$ with $Q_{2\rho}(z,\tau)\subset Q_{r}(x,t)$.

If $\epsilon_1=\epsilon_1(n,M_0)>0$ is small enough, then 
there exists $C=C(n,M_0,C_1,\zeta)$, 
$1\leq C<\infty$, $\gamma=\gamma(n, M_0, C_1, \zeta)>0$ and a 
Borel measurable function $g$ such that $d\mu /d\sigma=g$ on 
$\Sigma\cap Q_r(x,t)$ while
\begin{multline*}
\hspace{.1in}\left(\rho^{-n-1}\iint_{Q_{\rho}(z,\tau)\cap \Sigma} 
g^{1+\gamma}d\sigma\right)^{1/(1+\gamma)}\leq C\rho^{-n-1}
\iint_{Q_{2\rho}(z,\tau)\cap \Sigma }g\, 
d\sigma\\= C\rho^{-n-1}\mu\left(Q_{2\rho}(z,\tau)\right).
\end{multline*} 
\end{lemma}
\noindent Lemma \ref{lemma2} is a purely real variable result and therefore the proof can be readily adapted 
from the elliptic version given in \cite[Lemma 3.1]{BL}.   We omit the details.

\begin{proof}[Proof of Theorem \ref{tmain}]
Fix $Q_r(x_0,t_0)$ with $(x_0,t_0)\in\Sigma, r>0$, and let
$(Y,s)\in \Omega\setminus Q_{4r}(x_0,t_0)$. 
Then for all $Q_{\rho}(z,\tau)$ such that $Q_{2\rho}(z,\tau)\subset Q_r(x_0,t_0)$ we see that 
$(Y,s)\in \Omega\setminus Q_{4\rho}(z,\tau)$. 
Therefore Lemma \ref{lemma1} applies in each such
$Q_{\rho}(z,\tau)$ and if we set
$\mu:=\omega^{Y,s},$ then $\mu$ satisfies Lemma \ref{lemma2} with $C_1=C(\epsilon_1)$ 
and $\zeta=\eta(\epsilon_1)$ (here we are using
Lemma \ref{essbclosed}). Applying Lemma \ref{lemma2} we obtain Theorem \ref{tmain}. 
\end{proof}

 Before proceeding further, let us make the following geometric observation.  Set $R_0:=\diam(\Sigma)$.
Then there is a constant $c$,
depending only on dimension and ADR, such that $T_{max} - T_{min} \geq cR_0^2$.  Indeed, 
suppose first that $R_0<\infty$, and set
$r^2:=T_{max} - T_{min}$.  Then 
$\Sigma$ is contained in a closed rectangle in $\ree$ with dimensions $R_0\times...\times R_0 \times r^2$,
of volume $(R_0)^nr^2$.  We may then cover $\Sigma$ by a collection $\{Q^i_r\}_i$ of cardinality at most
$C (R_0/r)^n$, where for each $i$, $Q_r^i$ is a
parabolic cube of sidelength $r$, centered on $\Sigma$.  By the ADR property, since $\Sigma$ has diameter $R_0$,
we have
$$(R_0)^{n+1} \lesssim \sigma(\Sigma) \leq \sum_i \sigma(\Sigma \cap Q_r^i)
\lesssim \left(\frac{R_0}{r}\right)^n r^{n+1}\,.$$ 
Thus $r \gtrsim R_0$, as claimed.

Next, suppose that $R_0 = \infty$, but that $T_{max} -T_{min} =:r^2<\infty$. 
Then for any fixed $(x_0,t_0)\in\Sigma$, and any $R\in (r,\infty)$, the surface cube
$\Delta_R(x_0,t_0):= Q_R(x_0,t_0)\cap\Sigma$
is contained in a rectangle with dimensions $R\times...\times R\times r^2$, and volume $R^nr^2$.
We may then repeat the previous argument to see that $r\gtrsim R$, and then let $R\to \infty$. 


In the sequel, we shall continue to use the notation $R_0:= \diam(\Sigma)$. 

We now formulate the equivalence between the weak-$A_\infty$ property of parabolic measure,
and $L^p$ solvability of the initial-Dirichlet problem.

\begin{theorem}\label{t2.9}
Let $L$ be a divergence form parabolic operator defined on $\Omega$.
Suppose that  $\Sigma$ is globally time-backwards ADR, 
and assume further that if $R_0 = \infty$, then $T_{min} = -\infty$.
Then TFAE:


\begin{enumerate}
\item 
For every $\kappa_0 \in (0,1)$,  
 there is an exponent
$q>1$, possibly depending on $\kappa_0$,  such that  $\hm^{Y,s} \ll \sigma$, and
$k^{Y,s}:= d\hm^{Y,s}/d\sigma$ satisfies the reverse H\"older estimate
\begin{equation}\label{eq2.wRH}
\left(\, \,\tiltfiint_\Delta \left(k^{Y,s}\right)^q d\sigma \right)^{1/q} \lesssim_{\kappa_0} \,
\, \tiltfiint_{2\Delta} k^{Y,s} \,d\sigma\, 
\approx\,  \frac{\hm^{Y,s}(2\Delta)}{\sigma(\Delta)}\,,
\end{equation}
on every 
$\Delta =\Sigma \cap Q_r(x_0,t_0)$, with $t_0-T_{min} \geq \kappa_0 R_0^2$ and $r<\sqrt{\kappa_0} R_0/2$,
and for all
$(Y,s) \in \Omega\setminus Q_{4r}(x_0,t_0)$, 
uniformly for all such $\Delta$ and $(Y,s)$. 
\item   For every $\kappa_1 \in (0,1)$, there is an exponent $p<\infty$,  possibly depending on $\kappa_1$,
such that if 
$T_0-T_{min} \geq \kappa_1 R_0^2$,
and $f\in C_c(\Sigma)$ with compact support in $\Sigma^{T_0}$,  then
the parabolic measure solution 
$u$ of the initial-Dirichlet (resp., Dirichlet) 
problem  for
$L$ in $\Omega=\Omega^T$ with $T=T_{min}>-\infty$ (resp., in $\Omega$ if $T=-\infty$), 
with data $f$, satisfies 
for all $(x,t)\in \Sigma^{T_0}$, 
\begin{equation}\label{eq2.10}
N_* u(x,t)\, \lesssim_{\kappa_1} \,\left( \m(|f|^p)(x,t)\right)^{1/p}\,,
\end{equation}
where $\m$ denotes the parabolic Hardy-Littlewood maximal operator on $\Sigma$.
\item   For every $\kappa_1\in (0,1)$, there is an exponent $p<\infty$, possibly depending on
$\kappa_1$, such 
that if $T_0-T_{min} \geq \kappa_1 R_0^2$,
then the initial-Dirichlet problem for $L$ is $L^p$ solvable in $\Omega^{T_0}$,
with the estimate $\|N_*u\|_{L^p(\Sigma)} \lesssim_{\kappa_1} \|f\|_{L^p(\Sigma)}$. 
\end{enumerate} 
Furthermore, for appropriately related $\kappa_0, \kappa_1$, 
the exponents $q$ in (1), and $p$ in (2) and (3), satisfy the duality relationship $p=q/(q-1)$.
\end{theorem}

Here $N_* u$ denotes the non-tangential maximal function, of course taken with respect to
parabolic cones.  Precise definitions will be given in Section \ref{q}.

We note that we are implicitly assuming here, as above, that the continuous Dirichlet problem is solvable for
$L$;  we know that this is true if $L$ is the heat operator, or a divergence form parabolic operator
with $C^1$-Dini coefficients: see Remarks \ref{r1.22} and \ref{r-adr2}.

A few words of explanation are in order. 
 In less austere settings, say in Lipschitz cylinders or even
Lip(1,1/2) domains, the equivalence between $L^p$ solvability of the initial-Dirichlet problem and 
the $A_\infty$ property of
$\hm$, is well-known (see, e.g., \cite[Theorem 4.3]{N}:  for such a domain $\Omega$, with
$ T_{min}(\Omega)> -\infty$, one may consider the initial-Dirichlet problem 
in $\Omega = \Omega^T$, 
$T= T_{min}(\Omega)$, and then prove the main implication (1) implies (2) 
(or something essentially equivalent, namely that $N_*u \approx \m_{\hm} f$, where $\m_\hm$
 is the Hardy-Littlewood maximal operator with respect to parabolic measure at some fixed pole)
either by using Harnack's inequality and the
Harnack chain property to move from an arbitrary
point in a non-tangential ``cone" to a time forward point, or by
extending backwards in time, either by even reflection of the domain across the hyperplane 
$\{t\equiv T_{min}\}$, or simply by extending the 
time-slice $\Omega_{T_{min}}$ 
backwards in time as a cylinder.
In the more general setting considered in the present paper, neither of these technical
devices is available, and that is why we work in an ambient domain $\Omega$, and then 
solve the initial-Dirichlet problem in subdomains $\Omega^T$ with $T-T_{min} \gtrsim R_0^2$.

We mention that in Theorem \ref{t2.9}, we consider only the issue of $L^p$ {\em solvability}, 
i.e, existence of  solutions with $L^p$ estimates; we do not address the question of uniqueness.

\section{Proof of Lemma \ref{lemma1}}\label{s3}
\begin{proof}[Proof of Lemma \ref{lemma1}]
The proof is an adaption of the argument in \cite{BL} to the parabolic setting.
The principal new difficulty is the time lag inherent in parabolic problems.

Fix $\epsilon >0$, $r>0$, and $(x_0,t_0)\in\Sigma$, with
$0<r< \sqrt{t_0-T_{min}}/(8\!\sqrt{n})$,
 and suppose that
$\Sigma$ is TBADR on $\Delta_{2r}=\Delta_{2r}(x_0,t_0) = Q_{2r}(x_0,t_0) \cap\Sigma$.  
Observe that if \eqref{eq2.4} is true for 
some $\epsilon>0$, then it is true for all $\tilde{\epsilon}>\epsilon$.  Thus, 
we may suppose that $\epsilon \leq \epsilon_0$ for some sufficiently 
small but fixed
$\epsilon_0>0$.

Recall that $\delta({\bf Y}) := \dist({\bf Y},\eo)$, where the distance is of course the parabolic distance.
Replacing $Q_r$ by $Q_{2r}$ in Remark \ref{r1.eqdist}, we have
\begin{equation}\label{eqdist}
r< \sqrt{t_0-T_{min}}/(8\!\sqrt{n}) \,\implies\, 
\delta({\bf Y}) = \dist({\bf Y}, \Sigma)\,, \quad \forall\, {\bf Y} \in \Omega \cap Q_{2r}(x_0,t_0)\,.
\end{equation}
Thus, by hypothesis, we shall be working in a regime where $\delta({\bf Y}) = \dist({\bf Y}, \Sigma)$.
Moreover, we note that in this regime, i.e., for $(Y,s)\in \Omega \cap Q_{2r}(x_0,t_0)$, we have
\begin{equation}\label{eqdist2}
\dist\big((Y,s),\Sigma\big)=\delta(Y,s) \approx_a \dist\big((Y,s),\po\big)\,,\quad \text{if }\, s<t_0-(ar)^2\,,
\end{equation}
since $(x_0,t_0) \in \Sigma$ implies that $t_0\leq T_{max}$.  
Let us further note that $s<T_{max}$, 
for $(Y,s)\in \Omega \cap Q_{2r}(x_0,t_0)$, hence by Remark \ref{r-adr2},
\begin{equation}\label{eqdist3}
Q_\rho^-(Y,s) \cap \so =\emptyset\,,\quad 0<\rho<r/2\,,\,\,\, (Y,s)\in \Omega \cap Q_{3r/2}(x_0,t_0)\,.
\end{equation}


We shall use these facts repeatedly in the sequel, often implicitly.

 Let $M$ be a large positive constant to be chosen later.  Since
 it suffices to work with suitably small $\epsilon$, we may suppose that
$\epsilon\leq M^{-2}$.   Let $j$ be the greatest integer $\leq M/\epsilon$. Let
\begin{align}\nonumber
r_k^*&=\left(\dfrac{5}{4}+\dfrac{k}{4j}\right)r,
\\r_k'&=\left(\dfrac{5}{4}+\dfrac{k+1}{4j}\right)r,  \nonumber
\\\widehat{r}_k&=\left(\dfrac{5}{4}+\dfrac{k+1/2}{4j}\right)r. \label{eq3.1a}
\end{align}
Then define
\begin{align*}
U_k:= Q_{r_k'}(x_0,t_0)\setminus Q_{r_k^*}(x_0,t_0), \quad \text{ and} \quad
S_k:=\partial Q_{\widehat{r}_k}(x_0,t_0) \cap\Omega \,.
\end{align*}
for $1\leq k \leq j-1$. Note that $S_k\subset U_k$ and $U_k \subset Q_{\frac{3}{2} r}(x_0,t_0)$ for
each $k$.

Let $\epsilon' :=\dfrac{a}{MM_1}\epsilon$, and let $F\subset \Delta_{2r}$, with 
\[\sigma(F)\geq (1-\eta)\sigma ( \Delta_{2r})\,,\]
where $\eta =\eta(\epsilon)>0$ will be chosen momentarily.
We begin with two preliminary observations.

\begin{remark}\label{remark1}
First, assume that $(Y_1,s_1)\in Q_{\frac{7}{4}r}(x_0,t_0)$ with 
$$\dfrac{\epsilon'}{400}r\leq \delta(Y_1,s_1) \leq 
\dfrac{\epsilon '}{100}r\,.$$
Recalling \eqref{eq1.3} and Remark \ref{r2.5}, we set 
$Q_{Y_1,s_1} := Q\big((Y_1,s_1), 
20\frac{M_1}{a}\delta(Y_1,s_1)\big)$, so that $\Delta_{Y_1,s_1}=
\Sigma\cap 
Q_{Y_1,s_1}$, and note that $Q_{Y_1,s_1}\subset Q_{\frac{15}{8}r}(x_0,t_0)$. 
Then for $\eta(\epsilon)$ small enough,
\begin{align*}
\sigma\left(F\cap Q_{Y_1,s_1}\right)\geq (1-\theta)\sigma\left(\Delta_{Y_1,s_1} \right).
\end{align*}
Hence, by the local-ampleness assumption, setting $C=1/\beta$, we have
\begin{align}\label{eq3.2}
C\omega^{Y_1,s_1}(F)\geq 1,\quad (Y_1,s_1)\in Q_{\frac{7}{4}r}(x_0,t_0)\cap
 \{(Y,s):\dfrac{\epsilon'}{400}r\leq\delta(Y,s) \leq \dfrac{\epsilon '}{100}r\}\,.
\end{align} 
 \end{remark}
 
 \begin{remark}\label{remark2}
 Next, suppose that $(Y,s)\in S_k$, with
 $\delta(Y,s)<\epsilon'r/200$.
Then 
$$Q_{Y,s}:= Q((Y,s), 20\frac{M_1}{a}\delta(Y,s))\subset U_k\,.$$
 By Lemma \ref{Bourgain},
$$\omega^{Y,s}(Q_{Y,s} \cap \Sigma)\gtrsim 1\,,$$
and therefore, 
\begin{equation}\label{eq3.3}
\hspace{.1in}\omega^{Y,s}(U_k)\gtrsim 1\,,
\end{equation}
for $(Y,s) \in S_k\cap \{(Y,s):\, \delta(Y,s)<\epsilon'r/200\}$.
\end{remark}

We  now consider several cases.  Recall that
we have oriented our coordinate axes so that time runs from left to right.
Given a cube $Q=Q_\rho(x_0,t_0)$,
we shall use the terminology  
``back face of the boundary of $Q$" to denote the left-most face of $\partial Q$,
 i.e., the face with $t=t_0-\rho^2$.

\noindent \textbf{Case 1}. There is a point $(Y_0,s_0)$ on the back face 
of $S_k$ such that $\delta(Y_0,s_0)=\dfrac{\epsilon'}{200}r.$ 
Then
\begin{align}\label{eq3.3a}
\dfrac{\epsilon'}{400}r\leq \delta(Y_0, s_0-c_1r^2)\leq \dfrac{\epsilon'}{100}r,
\end{align}
where $c_1:=\left(\dfrac{\epsilon'}{200}\right)^4$, 
and $(Y_0,s_0-c_1r^2)\in  Q_{\frac{7}{4}r}(x_0,t_0)$. 
Hence, by \eqref{eq3.2}, 
\begin{align*}
C\omega^{Y_0,s_0-c_1r^2}(F)\geq 1.
\end{align*}

\smallskip

\noindent {\bf Claim 1}.  In the scenario of Case 1, 
given $\epsilon>0$, there exists a uniform constant $C_\epsilon$ such that
\begin{equation}\label{eq3.1}
C_{\epsilon}\omega^{Y,s}(F)\geq 1, \qquad \forall\,
(Y,s)\in S_k\cap \big\{(Y,s):\delta(Y,s) \geq \dfrac{\epsilon'}{200}r\big\} \,.
\end{equation}
\begin{proof}[Proof of Claim 1]
Let $(Y,s)\in S_k$ with $\delta(Y,s)\geq\dfrac{\epsilon'}{200}r$.

\noindent{\bf Case 1a}. $\delta(Y,s)\leq \dfrac{\epsilon'}{100}r.$ 
Then we are done by \eqref{eq3.2} (in this case with no dependence on $\epsilon$).

\noindent{\bf Case 1b}. $\delta(Y,s)> 
\dfrac{\epsilon'}{100}r$. 
 Note that $s\geq s_0$, since $(Y_0,s_0)$ lies on the back face of
$S_k$. Consider the 
parabola, call it $\C$, 
with vertex $(Y_0,s_0-c_1r^2)$, which opens to the right, 
contains the point $(Y,s)$, and has aperture $1/\alpha$, with
\begin{center}$\alpha=\dfrac{s-(s_0-c_1r^2)}{|Y_0-Y|^2}\geq 
\dfrac{c_1r^2}{|Y-Y_0|^2}\gtrsim c_1.$
\end{center}
Of course,
if $Y=Y_0$, then the parabola degenerates, and
$\C$ is  just a horizontal line, parallel to the $t$-axis. 
In any case, we start at the point $(Y,s)$ 
and move backward on $\C$, 
stopping the first time that we reach a point $(Z_0,\tau_0)$ with 
$\delta(Z_0,\tau_0)= \epsilon'r/100$. 
We eventually find such a point,  since 
$\delta(Y,s)>  \epsilon'r/100 
\geq \delta(Y_0,s_0-c_1r^2)$, by the scenario of Case 1b and \eqref{eq3.3a}. 
Note that by \eqref{eqdist3}, 
all the points on $\C$ between (and including) $(Y,s)$ and $(Z_0,\tau_0)$ lie in
$\Omega$, and thus in particular
\eqref{eq3.2} holds with $(Y_1,s_1) = (Z_0,\tau_0)$.
Moreover, since $||(Y,s)-(Z_0,\tau_0)||\lesssim r$, 
and $\delta(Z,\tau)\gtrsim_\epsilon r$  for every $(Z,\tau)\in\C$ between
$(Y,s)$ and $(Z_0,\tau_0)$, again using \eqref{eqdist3},
we can construct a Harnack path joining $(Z_0,\tau_0)$ to $(Y,s)$
along $\C$, to obtain \eqref{eq3.1} 
by Harnack's inequality
\cite[Theorem 1]{M}.  This proves Claim 1.
\end{proof}


\noindent {\bf Case 2}. For 
every point $(Y,s)$ on the back face of $S_k$, we 
have $\delta(Y,s)>\epsilon'r/200$.

In this case, we
slide the back face of $S_k$ forward in time
until we reach, for the first time, a point $(Y_0,s_0)$
with $\delta(Y_0,s_0)= \epsilon'r/200$.  Note that necessarily,  
\begin{equation}\label{eqso}
s_0  \leq \inf\left\{t: (x,t)\in \Delta_r(x_0,t_0) \text{ for some } x\right\} - (\epsilon'r/200)^2
\end{equation}
(otherwise, we would have stopped sooner when sliding the back face forward). 
In particular, by Remark \ref{1.8},
\[s_0 <  t_0-(ar)^2\,.\]
If we denote the boundary of the resulting rectangle 
by $S_k'$, then by construction
$(Y_0,s_0)$ is on the back face of $S_k'$.  
Since \eqref{eqdist2} holds for the point $(Y_0,s_0)$,
we may then repeat the argument
in Case 1, but with $S_k'$ in place of $S_k$.
 Consequently, we 
have the following.
\begin{remark}\label{remark3.5}
Estimate \eqref{eq3.1} holds, 
provided either that
$$(Y,s)\in S_k \cap \big\{(Y,s):\delta(Y,s) 
\geq \dfrac{\epsilon'}{200}r\big\}\,,\quad \text{or}\quad
(Y,s)\in S'_k \cap \big \{(Y,s):\delta(Y,s) 
\geq \dfrac{\epsilon'}{200}r\big\}\,,$$
in the scenarios of Case 1 and Case 2, respectively. 
Moreover, in Case 2,  $t=s_0$  on the back face
of $S_k'$, and by \eqref{eqso}, $s_0<t$ for any 
$(x,t) \in \Delta_r(x_0,t_0)=Q_r(x_0,t_0)\cap\Sigma$. 
\end{remark} 

\noindent {\bf Case 3}. For 
every $(Y,s)$ on the back face of $S_k$, we 
have $\delta(Y,s) < \epsilon'r/200.$ 

In turn, there are two sub-cases.

\noindent {\bf Case 3a}: For every point $(Y,s)$ on $S_k$, we 
have $\delta(Y,s)<\epsilon'r/200.$
In this case, the scenario of Remark \ref{remark2} applies to
every point $(Y,s)\in S_k$, and therefore \eqref{eq3.3} 
holds for all points on $S_k$.

\noindent {\bf Case 3b}. There exists a point $(Y^*,s^*)$ on $S_k$, with $\delta(Y^*,s^*)=\dfrac{\epsilon'}{200}r.$

Recall that $S_k = \partial Q_{\widehat{r}_k}(x_0,t_0)$, where
$\widehat{r}_k=\left(5/4+(k+1/2)/(4j)\right)r$. Consider the part of $Q_{\widehat{r}_k}(x_0,t_0)$ where 
$t_0-\widehat{r}_k^{\,\,2}<s<t_0-r^2$ and call this region $\widehat{I}_k$. 
Note that $|\widehat{I}_k|\approx r^{n+2}$. Cover $\widehat{I}_k$ by a union of 
pairwise non-overlapping half-closed
sub-cubes $\{Q^i\}_i$, such that $\epsilon'r\leq
\mathit{l}(Q^i)\leq 2\epsilon'r$.

\smallskip
\noindent {\bf Claim 2}: For $\epsilon$ small enough, at least one of the 
sub-cubes $Q^i$ misses $\Sigma$ (and thus also $\po$ by \eqref{eqdist} and the fact that
$t_0\leq T_{max}$).

Assume the claim momentarily.
Then there exists a point $(Z,\tau)\in Q^i$ such that 
$\delta(Z,\tau)>> \epsilon'r/200$. In the present scenario,
$\delta(Y,s)<\epsilon'r/200$ for every $(Y,s)$ on the 
back face of $S_k$, thus, there exist $\tilde{s}$ with
$s=t_0-\widehat{r}_k^{\,\,2}<\tilde{s}<\tau$, and $\delta(Z,\tilde{s})=
\epsilon'r/200.$ If we shrink $S_k$ by sliding the back face forward, 
increasing the time coordinate 
of the back face to $\tilde{s}$, and denote the boundary 
of the resulting rectangle by $S_k'$, 
then by construction, $(Z,\tilde{s})$ lies
on the back face of $S_k'$.  We can therefore
follow the argument given in 
Cases 1 and 2,
with $(Z,\tilde{s})$ playing the role of $(Y_0,s_0)$, to 
find that Remark \ref{remark3.5} 
continues to apply in Case 3b as well. 
We note in particular that by construction
$t= \tilde{s}$ on the back face of $S_k'$, and 
$\tilde{s} < t_0-r^2 <t$, for all $(x,t) \in Q_r(x_0,t_0)\cap\Sigma$. 

\begin{proof}[Proof of Claim 2] Suppose not. Then every $Q^i$ meets $\Sigma$. 
Hence there exists another parabolic cube $\widetilde{Q^i}$ 
such that $\widetilde{Q^i}\supset Q^i$, $\mathit{l}(\widetilde{Q^i})=4\mathit{l}(Q^i)$ and the 
center of $\widetilde{Q^i}$ is on $\Sigma$. By a 
rudimentary covering lemma argument, there 
exists a pairwise disjoint subcollection $\{\widetilde{Q}^{i_j}\}_j$,  with cardinality
$\#\{Q^{i_j}\}_j \approx \#\{Q^i\}_i$, such that
\begin{align*}
\bigcup_j 5\widetilde{Q}^{i_j} \supset \bigcup_i Q^i=\widehat{I}_k
\end{align*}
Let $\widetilde{I}_k$ be a fattened version of $\widehat{I}_k$, of 
comparable dimensions, such that $\bigcup_j\widetilde{Q}^{i_j} \subset \widetilde{I}_k$.
Then using disjointness and upper ADR, we obtain
\begin{equation}\label{eq3.4}
\sum_j \sigma(\widetilde{Q}^{i_j}\cap \Sigma)\leq \sigma(\widetilde{I}_k\cap \Sigma) 
\lesssim M_0 r^{n+1}\,.
\end{equation}
Using lower ADR, we obtain
\begin{equation}\label{eq3.5}
\sum_j \sigma(\widetilde{Q}^{i_j}\cap \Sigma) \gtrsim 
\#\{\widetilde{Q}^{i_j}\}\,(\epsilon' r)^{n+1},
\end{equation}
where as above, $\#\{\widetilde{Q}^{i_j}\}= $ cardinality of $\{Q^{i_j}\}.$
However,
\begin{align*}
\#\{\widetilde{Q}^{i_j}\}\approx 
\#\{Q^{i}\} \approx \dfrac{r^{n+2}}{(\epsilon'r)^{n+2}}\approx (\epsilon')^{-n-2}
\end{align*}
Therefore \eqref{eq3.5} becomes
\begin{equation}
\sum_j \sigma(\widetilde{Q}^{i_j}\cap \Sigma) 
\gtrsim (\epsilon')^{-n-2}(\epsilon'r)^{n+1}\approx \dfrac{r^{n+1}}{\epsilon'}>>M_0r^{n+1},
\end{equation}
for $\epsilon'\approx\epsilon$ small enough, contradicting \eqref{eq3.4}.
\end{proof}

Combining 
Remarks \ref{remark2} and \ref{remark3.5}, and our 
observation that the latter remark continues to hold in Case 3b,
we see that for $(Y,s) \in S_k$ (in Cases 1 and 3a), or  $(Y,s) \in S_k'$
(in Cases 2 and 3b),
\begin{equation*}
1\leq C\omega^{Y,s}(U_k)+C_{\epsilon}\omega^{Y,s}(F).
\end{equation*}
Moreover, letting $t_k$ denote the value of $t$ on the back face of $S_k$
or $S_k'$, as appropriate, we see that in every case,
$t_k<t$ for every $(x,t) \in \Delta_r(x_0,t_0):= Q_r(x_0,t_0)\cap\Sigma$.  
Consequently, by the weak maximum principle,
\begin{equation}\label{eq3.7}
\omega^{Y,s}\big(\Delta_r(x_0,t_0)\big)\leq C\omega^{Y,s}(U_k)+
C_{\epsilon}\omega^{Y,s}(F),
\end{equation}
for every $(Y,s) \in \Omega_k:= (\Omega \setminus R_k)\cap \{s>t_k\}$,
where $R_k$ is the closed cube, or rectangle, whose boundary is given by
$S_k$ or $S_k'$.
In addition, in every case, $t_k<t_0-(ar)^2$, and $R_k \subset Q_{4r}(x_0,t_0)$,
so in particular,  in \eqref{eq3.7},
we may take 
$$(Y,s)\in \big(\Omega\setminus Q_{4r}(x_0,t_0)\big) \cap \{s>t_0-(ar)^2\}\,.$$
For $(Y,s)$ in the latter set, we sum  \eqref{eq3.7} in $k$ to obtain
\begin{equation*}
\epsilon^{-1}
\omega^{Y,s}\big(\Delta_r(x_0,t_0)\big)\leq 
C\omega^{Y,s}(\Delta_{2r}(x_0,t_0))+
C_{\epsilon}\omega^{Y,s}(F),
\end{equation*}
since the sets $U_k$ are disjoint and are all contained in $Q_{2r}(x_0,t_0)$, and
 the cardinality of the index set $\{k\}$ is of the order $1/\epsilon$.
We now multiply by $\epsilon$ to reach the conclusion of the lemma in the special
case that $s>t_0-(ar)^2$.  

Let us now remove the latter restriction. 
Recall that the set $E(T)$ is defined in
\eqref{ETdef}. Observe that $t>t_0-r^2$ for every $(x,t) \in \Delta_r(x_0,t_0)$,
so that  if $s\leq t_0-r^2$, then
$\hm^{Y,s}\big(\Delta_r(x_0,t_0)\big) = 0$, and there is nothing to prove.  It therefore remains to treat the
case $t_0-r^2<s\leq t_0-(ar)^2$.    In this case, by an elementary covering argument,
we may cover the set $\Delta_r(x_0,t_0)\cap E(s)$ by a collection of surface cubes
$\{\Delta_i\}_{i=1}^N$, $\Delta_i = Q_i \cap \Sigma$, 
where $\Delta_i = \Delta_{cr}(x_i,t_i)$, with $(x_i,t_i)\in\Sigma$, $t_i\leq s$, and where
$c$ is a universal constant chosen small enough that $2Q_i\subset Q_{2r}(x_0,t_0)$; moreover,
this can be done in such a way that the cardinality $N$ of the collection is bounded by a universal
constant (depending on dimension).  
Thus, choosing $\eta'>0$ small enough, depending on $c$ and our previous choice of $\eta$, we
have that for $F\subset \Delta_{2r}(x_0,t_0)$,
\[\sigma(F) \geq (1-\eta')\, \sigma\big(\Delta_{2r}(x_0,t_0)\big) \implies
\sigma(F\cap2\Delta_i) \geq (1-\eta)\, \sigma\big(2\Delta_i\big)\,. \]
Since $s\geq t_i$ (hence in particular $s>t_i-(acr)^2$), 
we may therefore apply the previously treated special case to each $\Delta_i$, to
deduce that
\begin{multline*}
\hm^{Y,s}\big(\Delta_r(x_0,t_0)\big) = \hm^{Y,s}\big(\Delta_r(x_0,t_0)\cap E(s)\big)\\[4pt]
\leq \sum_{i=1}^N\hm^{Y,s}(\Delta_i)
\leq \epsilon \sum_{i=1}^N\hm^{Y,s}(2\Delta_i) +C_\epsilon  \sum_{i=1}^N\hm^{Y,s}(F\cap
2\Delta_i)\\[4pt]
\leq N \epsilon \hm^{Y,s}\big(\Delta_{2r}(x_0,t_0)\big) +NC_\epsilon \hm^{Y,s}(F)\,.
\end{multline*}
\end{proof}

\section{Proof of Theorem \ref{t2.9}}\label{q}

\begin{proof}[Proof of Theorem \ref{t2.9}]
Recall that either $L$ is the heat operator, or else we assume that
the continuous Dirichlet problem is solvable for $L$ in $\Omega$; in either case then,
the associated caloric/parabolic measure $\omega=\omega_L$ exists.  

The main step in the proof is to show that (1) implies (2).  
We turn our attention to this matter first.  The implication (2) implies (3) 
follows routinely from the density of $C_c(\Sigma)$ in
$L^p(\Sigma)$, and the self-improvement property of weak-reverse H\"older weights.  We omit the 
details, except to mention that in order for the non-tangential convergence to hold in a non-vacuous way,
one should impose some extra assumption to guarantee that the ``cones" defined in
\eqref{conedef} below are non-empty at infinitely many scales less than one,
$\sigma$ almost everywhere on $\Sigma$; an interior corkscrew 
condition is more than enough.
The implication (3) implies (1) will be proved at the end of this section.

\smallskip

\noindent {\em (1) implies (2)}.  Recall that $R_0:= \diam(\Sigma) \in (0,\infty]$.
We assume that for every $\kappa_0\in (0,1)$, and
for each $\Delta=Q_r(x_0,t_0)\cap\Sigma$, with $(x_0,t_0)\in\Sigma$,  
$t_0-T_{min} \geq \kappa_0 R_0^2$ and $r<\sqrt{\kappa_0} R_0/2$,
and for all 
$(Y,s)\in \Omega\setminus Q_{4r}(x_0,t_0)$, 
we  have  $\omega^{Y,s}\ll\sigma$ in $\Delta$,
and there exists some $q>1$ such that 
$k^{Y,s}\in$ weak-$RH_q(\Delta)$, with uniform constants, i.e. \eqref{eq2.wRH} holds for the Radon-Nykodym derivative
$k^{Y,s}:=d\omega^{Y,s}/d\sigma$, with $q$ and the implicit constants independent of $\Delta$ and $(Y,s)$;
equivalently, $\omega^{Y,s}\in$ weak-$A_{\infty}(\Delta)$, with uniform control of the constants.
Let $\kappa_1 \in (0,1)$, and set $\kappa_0:= \kappa_1/100$.
Suppose that $T_0-T_{min} \geq \kappa_1 R_0^2 =100 \kappa_0 R_0^2$.
Let $f$ be continuous on $\eo$,  with compact support in $\Sigma^{T_0}$,  and
 let $u$ be the parabolic measure solution 
 of the continuous Dirichlet problem (see Definition \ref{bvpdef}) for
$L$ with data $f$, in $\Omega$. 
For convenience, we shall treat the case $T_{min}>-\infty$; the proof in the case $T_{min} =-\infty$ 
is similar, but slightly simpler, and we leave  the details to the interested reader.

Our goal is to show that for $p=q/(q-1)$, and for all $(x,t)\in \Sigma^{T_0}$,
\begin{equation}\label{eq4.1}
N_* u(x,t) \lesssim\left( \m(|f|^p)(x,t)\right)^{1/p}\,,
\end{equation}
where $\m$ denotes the parabolic Hardy-Littlewood maximal operator on $\Sigma$, and hence that 
(1) implies (2).


To this end, we begin by defining non-tangential ``cones" and maximal functions, as follows. 
First, we fix a collection of 
 parabolic closed Whitney cubes covering $\Omega$, and we denote this collection by $\W$. 
 We also fix a constant
$\nu>0$ small enough so that, for each Whitney cube $I$,  its concentric parabolic dilate
$I^*:=(1+\nu)I$ 
also satisfies the Whitney properties. 
We will denote the collection of analogously
fattened Whitney cubes $\mathcal{W}^*$. 
 
 Given $(x,t)\in\Sigma$, set
\begin{equation}\label{wconedef}
\W(x,t):=\{I\in \W:  \diam(I) <10R_0 \,  \text{ and }\, \dist((x,t),I)\leq 100\diam(I)\},
\end{equation}
and define the (possibly disconnected) non-tangential ``cone" with vertex $(x,t)$ by
\begin{equation}\label{conedef}
\Upsilon(x,t):=\text{int}\left(\underset{I\in \W(x,t)}\bigcup I^*\right)\,,
\end{equation}
where int$(A)$ denotes the interior of the set $A$.
For a continuous $u$ defined on $\Omega$, the non-tangential maximal function of $u$ is defined by
\begin{equation}\label{ntdef}
N_*u(x,t):=\underset{(Y,s)\in \Upsilon(x,t)}\sup |u(Y,s)|.
\end{equation}

We now turn to the proof of \eqref{eq4.1}.


\noindent Splitting $f$ into its positive and negative parts, 
we may suppose without loss of generality that $f\geq 0$, hence also $u\geq 0$. 
Let $(x,t)\in\Sigma$ and fix $(Y_0,s_0)\in \Upsilon(x,t)$. Then $(Y_0,s_0)\in I_0^*$, 
for some $I_0\in \mathcal{W}(x,t)$ such that 
\begin{equation}\label{eq4.5aa}
r:= \delta(Y_0,s_0)\approx \diam(I_0)\approx ||(x,t)-(Y_0,s_0)||\,.
\end{equation} 
Of course, by definition of $\W(x,t)$ we have $r< K R_0$, for some sufficiently large universal constant $K$.

Let
$$Q_0:= Q\big((x,t), r\big)\,, \quad Q_k:= 2^kQ_0= Q\big((x,t), 2^kr\big)\,,\,\,\,\, k = 1,2,3...\,,  $$
and define corresponding surface cubes and subdomains:
\begin{equation}\label{okdef}
\Delta_k := Q_k\cap\Sigma\,,\quad \Omega_k:= Q_k\cap\Omega\,,\quad k \geq 0\,.
\end{equation}
Define a continuous partition of unity 
$\sum_{k\geq 0}\vp_k \equiv 1$ on $\Sigma$, such that $0\leq \vp_k\leq 1$ for all $k\geq 0$,
with
\begin{equation}\label{eq4.6}
\supp(\vp_0) \subset \Delta_2=:\RRR_0, \quad 
\supp(\vp_k) \subset \RRR_k:= \Delta_{k+2}\setminus \Delta_k,\,\, k\geq 1\,.
\end{equation}
 Set $f_k:= f\vp_k$, and let $u_k$ be the solution of the initial-Dirichlet problem in  $\Omega^T=\Omega$,
$T=T_{min}$, with initial data (at time $t=T_{min}$) equal to zero, and 
data $f_k$ on $\Sigma$.   Observe that since we are treating the case $T_{min}>-\infty$, and 
hence by assumption
$\diam(\Sigma)=R_0<\infty$, 
the boundary annulus $\RRR_k,\, k\geq 1$, is empty if $2^k r > R_0$;   for such $k$,
we have that $f_k$, and hence $u_k$, are identically zero.
Thus, we may restrict our attention to those $k$ for which
$2^kr\leq R_0$, so that $u =\sum_{  0\leq k\leq \log_2( R_0/r)} u_k$ in $\Omega $
(it may happen that $r>R_0$, but in that case only the term $k=0$ appears in the sum, and the 
following proof may be simplified considerably; we omit the very routine details).

Let us first observe that for each $k\geq 0$, and $(Y,s)\in \Omega$,
\begin{equation}\label{eq4.7}
u_k(Y,s) = \iint_{\RRR_k\cap\{T_0<\tau\}}f_k(z,\tau) \,d\hm^{Y,s}(z,\tau)\,.
\end{equation}
Indeed, since $f$ is supported in $\Sigma^{T_0}$, it follows that $u_k$ is zero
for $s\leq T_0$.

We now fix a sufficiently large integer $N$ to be chosen momentarily, and we
claim that
\begin{equation}\label{eq4.8}
\sum_{k=0}^N u_k(Y_0,s_0) \lesssim_N \left(\m\big(f^p\big)(x,t)\right)^{1/p}\,.
\end{equation}
To see this, we begin by recalling that by assumption, $T_0-T_{min} \geq 100 \kappa_0 R_0^2$, and
by construction, $r<KR_0$, for a suitably large universal constant $K$.
We may therefore cover 
$Q_{N+2} \cap\{(X,\tau):\,T_0<\tau\}$ by a collection $\F_0$ of pairwise disjoint half-open
parabolic cubes of parabolic sidelength $br$, 
 where $b$ is a sufficiently small number to be chosen momentarily, in particular
with $0<b\leq \sqrt{\kappa_0} /(100K)$.

Let $\F:= \{Q\in \F_0: Q \, \text{ meets } \, \Sigma\}$, and for each $Q\in \F$, let $\tQ$
be a cube centered at $(x_*,t_*)\in\Sigma$, 
containing $Q$,
of parabolic sidelength $\ell(\tQ) = 5\ell(Q) =5br$, so that $2Q_*\subset 100Q$.
Then for $b$ suitably small, $(Y_0,s_0) \in \Omega \setminus 4Q_*$ by \eqref{eq4.5aa}.
Thus, setting $\tD:= \tQ \cap\Sigma$, we note that by hypothesis, we may apply the reverse H\"older 
estimate \eqref{eq2.wRH} with $\Delta=\tD$, and with $(Y,s) = (Y_0,s_0)$.
Let $\tF$ denote the collection of all such $\tQ$.

Since $\sum_{k=0}^N f_k $ is supported in $\Delta_{N+2} = Q_{N+2}\cap \Sigma$, 
we then have that
\begin{multline*}
\sum_{k=0}^N u_k(Y_0,s_0)\, \leq\, \sum_{\tQ\in\tF} \iint_{\tD} f(z,\tau)
\,d\hm^{Y_0,s_0}(z,\tau)\\[4pt]
\leq \sum_{\tQ\in\tF} \sigma(\tD) \left(\,\, \tiltfiint_{\tD} f^p d\sigma\right)^{1/p}
 \left(\,\, \tiltfiint_{\tD} \left(k^{Y_0,s_0}\right)^q d\sigma\right)^{1/q} 
 \\[4pt] \lesssim
 \sum_{\tQ\in\tF}  \left( \,\,\tiltfiint_{\tD} f^p d\sigma\right)^{1/p},
\end{multline*}
where in the last step we have used the fact noted above: that the weak 
reverse H\"older estimate \eqref{eq2.wRH}
may be applied to each $\tD$ uniformly.  Now by the ADR property,  
$\sigma(\Delta_{N+3})\approx\sigma(\Delta_{N+2}) 
\approx_{N,n,K,\kappa_0} \sigma(\tD)$,
and by construction, we may suppose that each $\tD$ is contained in $\Delta_{N+3}$.  
Consequently, for each
$\tD$ we have
$$ \tiltfiint_{\tD} f^p d\sigma  \lesssim_{n,K,\kappa_0}
\tiltfiint_{\Delta_{N+3}}   f^p d\sigma
\lesssim_{n,K,\kappa_0} \m (f^p)(x,t)\,.$$
Since card$(\tF) \leq C(N,n,b)$, with $b$ in turn depending only on
$\kappa_0$ and $K$, the claimed bound \eqref{eq4.8} now follows.

Next, we claim that for $k\geq N+1$,  with $N$ chosen large enough, 
\begin{equation}\label{eq4.6a} u_k(Y_0,s_0) \lesssim \, 2^{-k\alpha}  \left( \m\big(f^p\big)(x,t)\right)^{1/p}\,.
\end{equation}
from which the desired bound \eqref{eq4.1} follows immediately, since $(Y_0,s_0)$ is an 
arbitrary point in $\Upsilon(x,t)$.

Recall that $\Omega_k$ is defined  in \eqref{okdef}.
We now fix $N$, depending only on the implicit constants in \eqref{eq4.5aa}, such that
$(Y_0,s_0) \in \Omega_{N-1}$.  Having fixed $N$,
we will allow implicit and generic constants to depend on $N$ without noting such dependence
explicitly.

For $k\geq N+1$, set
$$\W_k:= \big\{I\in\W:\, I\,\, {\rm meets}\,\, Q_{k-N}\big\}\,.
$$

Note that for $N$ chosen large enough, depending only on the Whitney construction,
we have that 
\begin{equation}\label{eq4.10}
\bigcup_{I\in \W_k} I^{*} \subset \Omega_{k-3} \,.
\end{equation}

Since $f_k$ vanishes in $\Delta_k$, by 
Lemma \ref{continuity}, we have
\begin{equation}\label{eq4.11}
u_k(Y_0,s_0) \lesssim 2^{-k \alpha } (2^kr)^{-n-2}  \iint_{\Omega_{k-N}}u_k \lesssim
 2^{-k \alpha } (2^kr)^{-n-2} \sum_{I\in\W_k} \iint_{I} u_k\,.
\end{equation}
Note that for each $I\in\W_k$, 
by the definition of $I^*$, we may fix a point $(Y_I,s_I) \in \partial I^*$ such that
$s_I > T_I + \nu \ell(I)^2$.  Note also that in particular, 
$(Y_I,s_I) \in \Omega_{k-3}$, by \eqref{eq4.10}.

For every $(Y,s) \in I$, by \eqref{eq4.7} we then have that by Harnack's inequality,
\begin{equation*}
u_k(Y,s) = \iint_{\RRR_k\cap\{T_0<\tau\}}f_k(z,\tau) \,d\hm^{Y,s}(z,\tau) 
\lesssim  \iint_{\RRR_k\cap\{T_0<\tau\}}f_k(z,\tau) \,d\hm^{Y_I,s_I}(z,\tau)\,.
\end{equation*}
Recall that since we are treating the case
$T_{min}>-\infty$, $R_0 <\infty$, we need only consider  $k$ such that 
$$2^kr \leq  R_0 \leq \sqrt{(T_0-T_{min})/(100\kappa_0)} \,.$$
We now choose a collection of surface cubes
$\F_k=\{\Delta^k_i= Q^k_i\cap\Sigma\}_i$, of parabolic sidelength $\ell(Q^k_i) = \sqrt{\kappa_0}\, 2^kr/100$, 
whose union 
covers $\RRR_k\cap\{T_0<\tau\}$. 
Note that we may do this in such a way that 
each $\Delta^k_i \subset \Delta_{k+3}\setminus \Delta_{k-1}$, and 
the cardinality of $\F_k$ is at most $C(n,\kappa_0)$, uniformly in $k$. 
Note further that by construction, the reverse H\"older estimate \eqref{eq2.wRH} may be applied uniformly 
in each $\Delta^k_i\in \F_k$, with pole at $(Y_I,s_I)$.
Consequently, for each $I\in \W_k$,
\begin{multline}\label{eq4.12}
\iint_I u_k \,\lesssim\, |I| \,\sum_{\F_k} \iint_{\Delta^k_i}f_k(z,\tau) \,d\hm^{Y_I,s_I}(z,\tau)\\[4pt]
\leq\,  |I|\, \sum_{\F_k} \sigma(\Delta^k_i) \left(\,\, \tiltfiint_{\Delta^k_i} f^p d\sigma\right)^{1/p}
 \left(\,\, \tiltfiint_{\Delta^k_i} \left(k^{Y_I,s_I}\right)^q d\sigma\right)^{1/q}\\[4pt] 
 \lesssim\, |I|\, \sum_{\F_k}  \left( \,\,\tiltfiint_{\Delta^k_i} f^p d\sigma\right)^{1/p}
 \,\lesssim_{\kappa_0} |I|\, \left( \m (f^p)(x,t)\right)^{1/p}\,,
\end{multline}
where in the last step we have used the bound on card$\left(\F_k\right)$, along with
ADR and the fact that $\diam(\Delta^k_i) \approx_{\kappa_0} \diam(\Delta_{k+3})$.
Note that \eqref{eq4.10} implies in particular that
$\sum_{I\in\W_k}|I| \lesssim (2^kr)^{n+2}$.
Plugging estimate \eqref{eq4.12} into \eqref{eq4.11},  we therefore
obtain the claimed bound \eqref{eq4.6a}, and hence that (1) implies (2).

\smallskip

\noindent {\em (3) implies (1)}.  We again treat only the case $T_{min} >-\infty$, 
$R_0<\infty$, as the proof in 
the case $T_{min} = -\infty$ is similar, but simpler.  Fix $\kappa_0 \in (0,1)$,  and a point
$(x_0,t_0)\in \Sigma$, with
$t_0 -T_{min} \geq \kappa_0 R_0^2$.  Let $0<r< \sqrt{\kappa_0}R_0/2$,
and set $\Delta=Q_r(x_0,t_0) \cap \Sigma$.   Our goal is to show that the reverse H\"older
estimate \eqref{eq2.wRH} holds for this $\Delta$, and to this end, it is actually enough to verify
\eqref{eq2.wRH} uniformly for each $\Delta'= Q_{\eps r}(x_1,t_1)\cap\Sigma \subset \Delta$, 
with $(x_1,t_1)\in \Delta$, where 
$\eps$ is a fixed small positive number to be chosen.  Indeed, the reverse H\"older
estimate for all such $\Delta'$,
with pole $(Y,s) \in \Omega\setminus Q_{4r}(x_0,t_0)$, 
implies that for $\Delta$, with constants depending on $\eps$. In turn, by Lemma \ref{lemma1}
and Lemma \ref{lemma2}, it is enough to show that there are uniform constants
$\theta,\beta \in (0,1)$ such that for every 
$(X,t) \in Q_{2\eps r}(x_1,t_1) \cap \Omega$,
if $E\subset \Delta_{X,t}$ is a
Borel set, 
\begin{equation}\label{eq4.13}
\sigma(E)\geq (1-\theta)\,\sigma(\Delta_{X,t}) \,\implies\,
\omega^{X,t}(E)\geq \beta\,,
\end{equation}
where as above, $\Delta_{X,t}:= \Sigma\cap Q_{X,t} :=\Sigma\cap Q((X,t), 20\frac{M_1}{a}\delta(X,t))$.

We fix $\Delta'$ and $(X,t)$ as above.
Let $(\hat{x},\hat{t})\in \Sigma$ be a touching point for $(X,t)$, i.e.,
$\delta(X,t) = \|(X,t) - (\hat{x},\hat{t})\|$.
We now choose $\eps$ small enough, depending on $M_1$, $a$, and $\kappa_0$, such that
for $(X,t) \in Q_{2\eps r}(x_1,t_1) \cap \Omega$, where  $(x_1,t_1)\in \Delta$ is the center of
$\Delta'$, we have
$200M_1 \delta(X,t)/a < \sqrt{\kappa_0}R_0/1000$, and also
$$\min (t-T_{min}, \hat{t} - T_{min}, T_{min}(\Delta')-T_{min}, T_{min}(\Delta_{X,t})-T_{min}) 
> \kappa_0 R^2_0/ 2\,.$$

Set $\kappa_1 = \kappa_0/100$, and set $T_0:= T_{min} + \kappa_1 R_0^2$.
Let $f \in C_c(\Delta_{X,t})$ be non-negative, with $\|f\|_{L^p(\Sigma)} \leq 1$.
By assumption, the solution $u$ to the initial-Dirichlet problem in
$\Omega^{T_0}$, with data $f$, enjoys the estimate
\begin{equation}\label{eq4.14}
\|N_*u\|_{L^p(\Sigma)} \lesssim_{\kappa_0} \|f\|_{L^p(\Sigma)} \approx_{\kappa_0} 1\,,
\end{equation}
for some $p<\infty$.

Set $r' := \delta(X,t)/10$.  Let $I\in \W$ be a Whitney cube containing $(X,t)$, and note that
$I\in \W(z,\tau)$, for every $(z,\tau) \in \Delta'':= Q_{r'}(\hat{x},\hat{t}) \cap\Sigma$ 
(see \eqref{wconedef}) and 
therefore
\begin{equation*}
u(Y,s) \lesssim\left(\,\, \tiltfiint_{\Delta''} \left(N_*u\right)^p d\sigma \right)^{1/p}\,,\quad \forall
(Y,s)\in I^*
\end{equation*}
(see \eqref{conedef}-\eqref{ntdef}).
Thus, by  \cite[Theorem 3]{M}, we have
\[u(X,t)\lesssim \left(\,\, \tiltfiint_{I^*} \big(u(Y,s)\big)^p dYds \right)^{1/p} 
\lesssim\left(\,\, \tiltfiint_{\Delta''} \left(N_*u\right)^p d\sigma \right)^{1/p}\lesssim
\delta(X,t)^{-(n+1)/p}\,,\] 
where in the last step we have applied the lower ADR estimate to $\Delta''$, and
used \eqref{eq4.14}. 
In turn, taking a supremum over all non-negative $f \in C_c(\Delta_{X,t})$ such that $\|f\|_p\leq 1$,
we obtain by Riesz representation that
\begin{equation}\label{eq4.15}
\left(\iint_{\Delta_{X,t}} \left(k^{X,t}\right)^q d\sigma \right)^{1/q} \lesssim \delta(X,t)^{-(n+1)/p}\,,
\end{equation}
with $q=p/(p-1)$.

We now claim that  the latter estimate implies  \eqref{eq4.13}, for 
suitable $\theta,\beta\in (0,1)$, in which case we are done. To prove this claim,
note first that by ADR, 
\begin{equation}\label{eq4.16}
\sigma(\Delta_{X,t}) \approx \delta(X,t)^{n+1}\,.
\end{equation} 
Let $E\subset \Delta_{X,t}$ satisfy the left hand estimate in \eqref{eq4.13}, and set
$A:= \Delta_{X,t} \setminus E$, for $\theta>0$ to be chosen,
so that 
\begin{equation}\label{eq4.17}
\sigma(A) \leq \theta \sigma(\Delta_{X,t})\,.
\end{equation}
Then
\begin{multline*}
\hm^{X,t}(A) \leq \sigma(A)^{1/p} \left(\iint_{\Delta_{X,t}} \left(k^{X,t}\right)^q d\sigma \right)^{1/q} 
 \\[4pt]
\lesssim \sigma(A)^{1/p}\delta(X,t)^{-(n+1)/p}
 \lesssim \theta^{1/p}
\approx \theta^{1/p} \hm^{X,t}(\Delta_{X,t})\,,
\end{multline*}
where in the last three steps we have used \eqref{eq4.15}-\eqref{eq4.17}, and then Lemma
\ref{Bourgain}.  Taking complements, and using Lemma \ref{Bourgain} once again,
for $\theta$ small enough we obtain \eqref{eq4.13}. 
\end{proof}

\section{Two Applications}\label{s5a}

In this section, we briefly discuss two applications of our main result, Theorem \ref{tmain}.
The first is an extension of a result proved in \cite{NS}.  We refer to \cite{NS} for definitions of the terms
not previously defined in the present paper.

\begin{theorem}\label{NSweak} Let $\Omega \subset \ree$ satisfy a time-synchronized two cubes condition, 
and suppose that $\Sigma$ is parabolic uniformly rectifiable (in particular $\Sigma$ is globally ADR).
Then caloric measure $\hm$ satisfies a local weak-$A_\infty$ condition with respect to $\sigma$, equivalently,
$\hm \ll\sigma$ and the Radon-Nikodym derivative $d\hm/d\sigma$ verifies the weak Reverse H\"older condition
\eqref{weakRH}.
\end{theorem}

A few remarks are in order.  In \cite{NS}, the authors obtain a similar result, but assuming in addition
that $\Omega$ satisfies a parabolic Harnack Chain condition.    Our Theorem \ref{tmain} allows us to dispense
with the latter connectivity assumption in Theorem \ref{NSweak}.  The conclusion in \cite{NS} is that
$\hm$ satisfies an $A_\infty$ condition (which entails doubling) with respect to $\sigma$, 
but in the absence of connectivity the non-doubling
weak-$A_\infty$/weak Reverse H\"older conclusion is best possible.
We remark that the aforementioned result of \cite{NS}, and our Theorem \ref{NSweak}, are the parabolic analogues of
results proved in \cite{DJ} and in \cite{BL}, respectively.

  The hypotheses of the theorem correspond to
the case $T_{min}=-\infty$, $T_{max} = \infty$.  
A sketch of the proof is as follows. One first invokes the deep fact proved in
\cite[Theorem 1.2]{NS} that under the hypotheses of Theorem \ref{NSweak}, one obtains an interior 
``big pieces" approximation (see \cite{NS} for the precise definition), analogous to that proved in the elliptic setting in  \cite{DJ}, by domains of the sort
considered in \cite{LM}. By the result of \cite{LM}, plus a standard maximum principle argument,
one obtains the $(\theta,\beta)$-local ampleness condition \ref{eq1.4}.  In addition, it is not difficult to show
that in the presence of ADR, the time-synchronized two cubes condition of \cite{NS} implies global time-backwards ADR
(in fact, it implies a time-symmetric version of ADR, in which one has thickness both in
$\Delta_r^-$ and in $\Delta_r^+$).  At this point, the conclusion follows by Theorem \ref{tmain}.

A second application of Theorem \ref{tmain}
will appear in our forthcoming joint paper \cite{GH}.  We state the result here, but refer the reader to that 
paper for details.

\begin{theorem} Let $\Omega\subset \ree$, whose quasi-lateral boundary $\Sigma$ is
globally ADR and time-backwards ADR.  Suppose that the Dirichlet problem is ``BMO-solvable" in $\Omega$.
Then caloric measure $\hm$ satisfies a local weak-$A_\infty$ condition with respect to $\sigma$, equivalently,
$\hm \ll\sigma$ and the Radon-Nikodym derivative $d\hm/d\sigma$ verifies the weak Reverse H\"older condition
\eqref{weakRH}.

\end{theorem}

The result is a parabolic version of the main theorem of \cite{HLe}, which in turn entailed removing
all connectivity hypotheses (in particular, the Harnack Chain condition), from earlier elliptic results of
\cite{DKP} and \cite{Z}.  We refer to \cite{GH} for a precise formulation of the BMO-solvability statement
in the parabolic setting.

\appendix
\section{Proof of Bourgain-type Estimate,  Lemma \ref{Bourgain}}\label{s4}

\begin{proof}
Let $(x_0,t_0)\in \Sigma$, 
and let $0<Mr< \sqrt{t_0-T_{min}}/(4\!\sqrt{n})$, where our goal is now,
equivalently, to show that 
$\omega^{X,t}\big(\Delta_{Mr}\cap E(T)\big)\geq \eta$,
for all $(X,t)\in Q_{ar/2}(x_0,t_0)\cap\Omega$,  
with $\Sigma$ time-backwards ADR on $\Delta_{Mr}=Q_{Mr}(x_0,t_0)\cap \Sigma$,
and $M$ is a large enough constant depending only on  $\lambda$, $n$ and 
ADR (including time-backwards ADR).  Here, $E(T)$ is defined as in
\eqref{ETdef}, with $T=T_{max}\big(Q_{ar/2}(x_0,t_0)\big)= t_0 +(ar)^2/4$.
We then obtain the conclusion of Lemma \ref{Bourgain}  with $M_1=2M$.

We continue to
assume either that
$L$ is the heat operator, or that the continuous Dirichlet problem is solvable in $\Omega$ for $L$.

\smallskip
\noindent {\bf Claim 1}.  There exist numbers  $a \in (0,1/2)$ and $b\in (0,1)$, depending only on $n$ and  ADR  (including
the time-backwards ADR constants), such that
if $\Sigma$ is time-backwards ADR on $ \Delta=\Delta_r =\Sigma \cap Q_r(x_0,t_0)$, then 
\begin{equation}\label{a1}
\sigma\left(Q^-_r\big(x_0,t_0-(ar)^2\big)\cap\Sigma\right)\geq 
\sigma\left(\Delta_r^-\cap \left\{t<t_0-(ar)^2\right\}\right) \geq br^{n+1}\,.
\end{equation} 
\begin{proof}[Proof of Claim 1] Observe that
$$Q^-_r\big(x_0,t_0-(ar)^2\big) \supset Q_r^-\cap \left\{t<t_0-(ar)^2\right\}\,.$$  
Thus,
the first inequality in \eqref{a1} is trivial, so we need only prove the second.
Set $\Phi_{ar}:=Q_r^-(x_0,t_0)\cap \{t_0-(ar)^2\leq t\leq t_0\}$ and take $a=2^{-m}$, where 
$m> 1$ will be chosen large enough.
We then decompose $\Phi_{ar}$ into sub-cubes, with parabolic length $ar$, all of equal size, and denote these cubes by $Q_{ar}^i$. As the $n$-dimensional measure of the face of each $Q_{ar}^i$ with $t=t_0$ is $(ar)^n$, and the $n$-dimensional measure of the face of $Q_r^-(x_0,t_0)$, again with $t=t_0$, 
is $r^n$, there are $a^{-n}$ such sub-cubes $Q_{ar}^i$. 
Hence,
\begin{align*}
\Phi_{ar}=\bigcup_{i=1}^{a^{-n}} Q_{ar}^i.
\end{align*}
By upper ADR of $Q_{ar}^i\cap\Sigma$, this yields in turn that
\begin{align*}
\sigma\left(\Sigma\cap \Phi_{ar}\right) 
\leq M_0\sum_{i=1}^{a^{-n}} (ar)^{n+1}\leq M_0ar^{n+1}.
\end{align*}
Therefore, if we choose $a\leq (2M_0)^{-1} c$, where $c$ is the constant in the definition of time-backwards ADR,
\begin{equation}\label{eq4.2}
\sigma\left(\Sigma\cap \Phi_{ar}\right) 
\leq \dfrac{1}{2}cr^{n+1}
\leq \dfrac{1}{2}\sigma\left(\Sigma\cap Q_r^-(x_0,t_0) \right).
\end{equation}
Consequently, since $Q_r^-(x_0,t_0)\setminus \Phi_{ar} 
= Q_r^- \cap \left\{t<t_0-(ar)^2\right\}$, 
\begin{align*}
\sigma\left(\Delta_r^- \cap \left\{t<t_0-(ar)^2\right\}\right)
&=
\sigma\left(\Sigma\cap \left(Q_r^-(x_0,t_0)\setminus \Phi_{ar}\right)\right)
\\&\geq \dfrac{1}{2}\sigma\left(\Sigma\cap Q_r^-(x_0,t_0)\right)
\geq \dfrac{1}{2}cr^{n+1},
\end{align*}
where we have used time-backwards ADR in the last step.
\end{proof}
To prove the estimate of Bourgain-type we will  use \eqref{a1}.  
We follow the proof in \cite{B}, adapting it to the parabolic setting.
Set
\begin{equation}\label{a3a}
\widehat{\Delta}:= Q_{r}^-(x_0,t_0-(ar)^2)\cap\Sigma\,,
\end{equation}
and define 
\begin{equation}\label{udef} u(X,t):=
\int_{-\infty}^t\int_{\Sigma_s} 
\Gamma(X,t,y,s)\,\chi_{\widehat{\Delta}}(y,s)\,d\sigma_s(y)ds\,,
\end{equation}
where $\Gamma(X,t,Y,s)$ is the fundamental solution of $L$. 
In \cite{QX}, the authors prove the following inequality:
\begin{align}\nonumber
\dfrac{1}{(N(t-s))^{n/2}}\exp\left(-\dfrac{N|X-y|^2}{t-s}\right)\chi_{\{t>s\}}\, &\leq\, \Gamma(X,t,y,s) 
\\
\label{a2}
& \leq\,\dfrac{N}{(t-s)^{n/2}}\exp\left(-\dfrac{|X-y|^2}{N(t-s)}\right)\chi_{\{t>s\}},
\end{align}
where $N$ depends on dimension and $\lambda$. 
\begin{remark} In fact, in \cite{QX}, 
the authors obtain this inequality in the more general situation that 
the non-symmetric part of the coefficient matrix $A$ belongs to $BMO$.
\end{remark}


\begin{claim}\label{c7} The function
$u$ defined in \eqref{udef} satisfies the following properties: there exist
constants $C_1$ 
and $c_2$ depending only on ADR, $\lambda$, and $n$, such that
\begin{enumerate}
\item $u$ is continuous in $\ree$, $Lu=0$ in $\Omega,$ and $u\equiv 0$ for $\{t: t\leq t_0-\frac{5}{4}r^2\}$.
\smallskip
\item $0\leq u \leq C_1 r$ in $\ree$. 
\smallskip
\item $u(X,t)\geq c_2 r$, for $(X,t)\in 
Q_{ar/2}\cap\overline{\Omega}$, where $Q_{ar/2}:= Q_{ar/2}(x_0,t_0).$
\smallskip
\item $u(X,t)\leq C_1 M^{-n}r$, for 
$(X,t)\in \overline{\Omega}\setminus Q_{Mr/2}$, 
where $Q_{Mr/2} := Q\left((x_0,t_0),\frac{M}2r\right)$.
\end{enumerate}
\end{claim}

Property (1) follows by definition of $u$, and the fact that $(ar)^2 \leq r^2/4$.  
Let us now verify the remaining properties.
 
\begin{proof}[Proof of property (2)]
Certainly, $u\geq 0$, by definition.  To prove the upper bound for $u$, we first note that
\begin{equation}\label{a4}
\Gamma(X,t,y,s)\lesssim_{n,\lambda} ||X-y,t-s||^{-n}\,,
\end{equation}
as one may readily derive from \eqref{a2} and the definition of the parabolic distance; we omit the details.

Next, we split $u$ as follows:
\begin{align*}
u(X,t)&=\displaystyle\iint_{\widehat{\Delta}}\Gamma(X,t,y,s)\,d\sigma(y,s)
\\&=\iint_{\widehat{\Delta}\cap \{||X-y,t-s||< r\}} \Gamma\,d\sigma
\,+\, \iint_{\widehat{\Delta} \cap \{||X-y,t-s||\geq r\}} \Gamma\,d\sigma\,
=:\,I+II.
\end{align*}
By \eqref{a4}, the integrand in term $II$ is at most $Cr^{-n}$, hence,
\begin{equation*}
II \lesssim 
r^{-n}\sigma(\widehat{\Delta})\lesssim r^{-n}r^{n+1}\leq C r,
\end{equation*}
by upper ADR,  where $C=C(\text{ADR},n, \lambda)$.

In term $I$,
we will dyadically decompose $\widetilde{\Delta}:=\widehat{\Delta}\cap \{||X-y,t-s||\lesssim r\}$. 
Define 
\[A_k:=\{(y,s)\in \widetilde{\Delta}: 2^{-k-1}r < ||X-y,t-s||\leq 2^{-k}r\}\,.\]
As $A_k\subset \widehat{\Delta}$, we can use upper ADR on $A_k$, along with
\eqref{a4}, to obtain the following estimate:
\begin{align*}
I&\lesssim \sum_{k=0}^{\infty}\iint_{A_k}||X-y,t-s||^{-n}d\sigma(y,s)
\\&\leq \sum_{k=0}^{\infty}(2^{-k-1}r)^{-n}\iint_{A_k}d\sigma(y,s)
\\&\lesssim \sum_{k=0}^{\infty} (2^{-k-1}r)^{-n}(2^{-k}r)^{n+1} 
\leq C r,
\end{align*}
 where $C=C(\text{ADR},n,\lambda)$. 
\end{proof}

\begin{proof}[Proof of property (3)]

As $(y,s)\in \widehat{\Delta}$ and $(X,t)\in Q_{\frac{a}{2}r}\cap \overline{\Omega}$, we have that 
\begin{center}$t_0-(ar)^2-r^2<s<t_0-(ar)^2$,
\end{center}
\begin{center}
$t_0-\dfrac{1}{4}(ar)^2<t<t_0+\dfrac{1}{4}(ar)^2.$ 
\end{center}
Hence
\begin{align*}
t-s&<t_0+\dfrac{1}{4}(ar)^2-(t_0-(ar)^2-r^2)=\dfrac{5}{4}(ar)^2+r^2, \quad \text{and}
\\t-s&>t_0-\dfrac{1}{4}(ar)^2-(t_0-(ar)^2)=\dfrac{3}{4}(ar)^2.
\end{align*} 
Overall, this gives us that $t-s\approx r^2,$
with implicit constants depending on $a$. Since $|X-x_0|\lesssim r$ and 
$|y-x_0|\lesssim r$, 
by the triangle inequality we have
$|X-y|\lesssim r$, hence, 
\[-|X-y|^2 \gtrsim -r^2.\] 
Combining the above estimates, as well as using ADR, we obtain that
\begin{align*}
u(X,t)&\geq\iint_{\widehat{\Delta}}\dfrac{1}{(N(t-s))^{n/2}}\exp\left(-\dfrac{N|X-y|^2}{(t-s)}\right)\chi_{\{t>s\}}
d\sigma(y,s)
\\&\gtrsim \iint_{\widehat{\Delta}} \dfrac{1}{r^n}\exp\left(-\dfrac{r^2}{Cr^2}\right)d\sigma(y,s)
\\&\gtrsim r^{-n}\sigma(\widehat{\Delta})\gtrsim r^{-n}r^{n+1}\geq c_2r,
\end{align*}
by \eqref{a1} and the definition of $\widehat{\Delta}$ \eqref{a3a}, where $c_2=c_2(n,ADR, \lambda)>0$
(recall that by Claim 1, $a$ depends only on $n$, ADR and time-backwards ADR). 
\end{proof}

\begin{proof}[Proof of property (4)]

As $(y,s)\in \widehat{\Delta}$ and $(X,t)\in \overline{\Omega}\setminus Q_{Mr/2}$, we have that
\begin{equation*}
\|(X,t)-(y,s)\|\geq \big|\|(X,t)-(x_0,t_0)\|-\|(x_0,t_0)-(y,s)\| \big|
\geq |cMr-Cr|
\, \gtrsim Mr,
\end{equation*}
provided that $M$ is chosen sufficiently large.
Combining the latter estimate with \eqref{a4}, we have
\begin{equation*}
u(X,t) 
\lesssim_{n,\lambda} (Mr)^{-n}\sigma(\widehat{\Delta}) 
\leq C_1 M^{-n} r,
\end{equation*}
where $C_1=C_1(n,ADR,\lambda)$.
\end{proof}

\begin{claim}\label{c9}
Set $\tilde{u}(X,t):=\dfrac{1}{r}\left(u(X,t)-\underset{\overline{\Omega}\setminus Q_{Mr/2}}\sup\, u \right).$ Then $\tilde{u}$ satisfies the following:
\begin{enumerate}[label=(\roman*)]
\item $\tilde{u}$ is continuous in $\ree$, $L\tilde{u}=0$ in $\Omega$, 
and $\tilde{u} \leq 0$ in $\overline{\Omega}\setminus Q_{Mr/2}$, and in $ \{t: t\leq t_0-\frac{5}{4}r^2\}$.
\smallskip
\item $|\tilde{u}(X,t)|\leq 2 C_1$ in $\overline{\Omega}$.  
\smallskip
\item $\tilde{u}(X,t)\geq \frac{1}{2}c_2$ for $(X,t)\in Q_{ar/2}\cap\overline{\Omega}$, provided $M$ is large enough. 
\end{enumerate}
\end{claim}
Indeed, 
property (i) follows immediately from Claim \ref{c7}, property $(1)$, and
the definition of $\tilde{u}$;
property (ii)  from
Claim \ref{c7}, property (2), and property (iii) from Claim \ref{c7}, properties (3) and (4).  
We omit the routine details.

Recall that $E(T)$ is defined as in \eqref{ETdef}, with $T:=T_{max}(Q_{ar/2}) = t_0+(ar)^2/4$.


\smallskip
\noindent {\bf Claim 2}.  
$\omega^{X,t}\big(\Delta_{Mr}\cap E(T)\big)\gtrsim \tilde{u}(X,t)$, for
\[(X,t) \in  \widetilde{\Omega}:= \Omega \cap Q_{Mr/2} \cap  \{t: t> t_0-5r^2/4\}\cap E(T)\,.\]
\begin{proof}[Proof of Claim 2]
Let $0<\eps\ll r$.  Set $\widetilde{\Omega}_\eps:=  \widetilde{\Omega} \cap E(T-\eps)$.
Observe that by property (i) and the definition of $\widetilde{\Omega}_\eps$,
\[\mathcal{P}\widetilde{\Omega}_\eps \subset \{  \tilde{u} \leq 0\} 
\cup \left(\overline{Q_{Mr/2} }\cap\Sigma\cap E(T)\right)
\subset  \{  \tilde{u} \leq 0\} \cup \left(\Delta_{Mr}\cap E(T)\right)\,.\]
The claim then follows with $\widetilde{\Omega}_\eps$ in place of
$\widetilde{\Omega}$, by property (ii) and the weak maximum principle.  The full
claim follows by letting $\eps \to 0$,
\end{proof}

Note that $Q_{ar/2}\cap\Omega \subset \widetilde{\Omega}$, since $t_0 -(ar)^2/4<t<T$ in
$Q_{\frac{a}{2}r}$.
By  Claim \ref{c9}, property (iii), we have that $\tilde{u}(X,t)\geq \frac{1}{2}c_2$ for 
$(X,t)\in Q_{ar/2}\cap\Omega$. 
Thus, for such $(X,t)$, by Claim 2 we obtain  
\begin{align*}
\frac{1}{2}c_2 \leq \tilde{u}(X,t) &\lesssim \omega^{X,t}(\Delta_{Mr}).
\end{align*}
This finishes the proof with $\eta \approx c_2.$
\end{proof}

\section{Proof of H{\"o}lder Continuity at the Boundary, Lemma \ref{continuity}}\label{s5}

As above, given a cube $Q_r$ centered on $\Sigma$, and a fixed time $T$,
we set $\Omega_r:= Q_r\cap\Omega$,
and $\Omega_r(T):= \Omega_r\cap E(T)$, where we recall that $E(T)$ is defined in \eqref{ETdef}.
We first state a version of Bourgain's lemma for supersolutions.
\begin{lemma}[Parabolic Bourgain-type Estimate for supersolutions]\label{Bourgainsup} 
Let $(x_0,t_0)\in\Sigma$, and let $0<r< \sqrt{t_0-T_{min}}/(4\!\sqrt{n})$. 
Set $Q_r:= Q_r(x_0,t_0)$, $Q_{ar/M_1}=
Q((x_0,t_0),\frac{a}{M_1}r)$, and define
\[T=T_{max}(Q_{ar/M_1})=
t_0+(aM_1^{-1}r)^2\,.\]
Assume that 
$\Sigma$ is time-backwards ADR on $ \Delta_r := Q_r \cap \Sigma$.
Then there exists $M_1,\eta>0$ such that 
if $w$ is a  non-negative
supersolution in $\Omega_r(T)$, with $w\geq 1$ on $\Delta_r\cap E(T)$
in the sense that
 \[\liminf_{(X,t)\to (y,s)} w(X,t) \geq 1\,,\quad (y,s) \in \Delta_r\cap E(T)\,,\]
then 
\begin{equation*}
w(X,t)\geq \eta\,,\quad \forall \, (X,t)\in \Omega_{ar/M_1}\,.
\end{equation*}
\end{lemma}

The proof is identical to that of Lemma \ref{Bourgain}:  replace $r$ by $M r$ (with $2M=M_1$), 
build the same auxiliary solutions
$u$ and $\tilde{u}$, and then apply the weak maximum principle as before.  We omit the details.

We next establish the following H\"older continuity statement for subsolutions.

\begin{lemma}\label{lb2}
Suppose that 
$\Sigma$ is time-backwards ADR on $Q_{2r}:=Q_{2r}(x_0,t_0)$, with $(x_0,t_0)\in \Sigma$ and 
$0<r <\sqrt{t_0-T_{min}}/(8\!\sqrt{n})$.   
Then there exists $C=C(n,\text{ADR}), \alpha=\alpha(n,\text{ADR})>0$, 
such that  if $v$ is a
non-negative subsolution in $\Omega_{2r}(T_1)$, which 
vanishes continuously on
 $\Delta_{2r}\cap E(T_1):=Q_{2r}\cap \Sigma \cap E(T_1)$, then 
\begin{equation*}
v(Y,t)\leq C\left(\dfrac{\delta(Y,t)}{r}\right)^{\alpha} \m(v) \,,
\quad \forall \,(Y,t)\in \Omega_{r} \,, 
\end{equation*}
where 
$T_1:= T_{max}(Q_r) = t_0+r^2$, and $\m(v):= \sup_{\Omega_{3r/2}(T_1)} v$.
\end{lemma}

\begin{proof}[Proof of Lemma \ref{lb2}] If  $\m(v) =\infty$,  there is nothing to prove, 
so we may assume that $\m(v) < \infty$. Normalize $v$ 
so that $\m(v) \leq 1$, and set $w:=1-v$. Then 
$0\leq w\leq 1$ in $\Omega_{3r/2}(T_1)$, and
$w\equiv 1$ on $\Delta_{2r}\cap E(T_1)$.
Then by Lemma \ref{Bourgainsup}, there exists $M, \eta >0$ such that
$w(Y,t)\geq \eta,$
and therefore
\begin{align*}
v(Y,t)\leq 1- \eta\,,\quad \forall\, (Y,t)\in Q_{\frac{a}{M}r}(x_0,t_0)\cap \Omega\,.
\end{align*}
Iterating, 
and using the fact that similar estimates hold with
$Q_{3r/2}(x_0,t_0)$ replaced by $Q_{r/2}(x_1,t_1)$, for $(x_1,t_1)\in \Delta_{r}(x_0,t_0)$, we obtain the 
conclusion of the lemma.
\end{proof}

\begin{proof}[Proof of Lemma \ref{continuity} (H{\"o}lder Continuity at the Boundary)]
Set 
\[T_1 := T_{max}(Q_r(x_0,t_0) = t_0+r^2\,,\]
Set $\m(u):= \sup_{\Omega_{3r/2}(T_1)} u$.
The first step is to establish the estimate
\begin{equation}\label{eqb3}
u(Y,t)\leq C\left(\dfrac{\delta(Y,t)}{r}\right)^{\alpha}\m(u)\,,
\quad \forall (Y,t)\in \Omega_r\,,
\end{equation}
where $u$ is the parabolic measure solution with non-negative data $f\in C_c(\eo)$, with
$f\equiv 0$ on $\Delta_{2r}$.  
If $L$ is a divergence form parabolic operator for which
the continuous Dirichlet problem is solvable, then  estimate \eqref{eqb3} is a special case of
Lemma \ref{lb2}, since a non-negative solution is in particular a non-negative subsolution.

On the other hand, suppose now that $L$ is the heat operator. In this case,
we  need not assume {\em a priori}
solvability of the continuous Dirichlet problem; rather
we shall use Lemma \ref{lb2}, and the Perron construction (see \cite[Chapter 8]{W2}).  By resolutivity of
$C(\eo)$, the caloric measure solution is the Perron (more precisely, PWB) solution, and is given by
\[u = \sup \left\{v:\,v\in\mathcal{L}_f\right\},\]
where, since $f \geq 0$, without loss of generality the lower class $\mathcal{L}_f$ consists of all 
{\em non-negative} subcaloric
$v$ satisfying
 \[\limsup_{(X,t)\to (y,s)} v(X,t) \leq f(y,s)\,,\quad (y,s) \in \no\,,\]
and
 \[\limsup_{(X,t)\to (y,s^+)} v(X,t) \leq f(y,s)\,,\quad (y,s) \in \sso\,.\]
 In particular, since $\Delta_{2r} \subset \Sigma \subset \no$ (recall that $\sso \cap \Delta_{2r} = \emptyset$,
 by the time-backwards ADR assumption; see Remark \ref{r-adr2}), we have
  \[0 \leq \limsup_{(X,t)\to (y,s)} v(X,t) \leq f(y,s) = 0 \,,\quad (y,s) \in \Delta_{2r}\,,\]
  for all $v \in \mathcal{L}_f$.  
  Thus, each such $v$ vanishes continuously  on $\Delta_{2r}$, and therefore Lemma \ref{lb2} may be applied to
  any $v\in \mathcal{L}_f$.  
  
  If $\m(u)=\infty,$ then \eqref{eqb3} is trivial, so after normalizing,
we may suppose that $\m(u)\leq 1$. By the PWB construction, we then have
$\m(v):=\sup_{\Omega_{3r/2}(T_1)} v\leq 1$, for all $v \in \mathcal{L}_f$.
  Given $\eps>0$, and a point $(Y,t)\in  \Omega_r$, 
  we may choose $v\in \mathcal{L}_f$ such that
  \[ u(Y,t)  \leq v(Y,t) +\eps\,.\]
  Applying Lemma \ref{lb2} to $v$, with $\m(v)\leq 1$,
  and letting $\eps\to 0$, we obtain \eqref{eqb3}.
  
  With \eqref{eqb3} in hand, it remains to replace $\m(u)$ by an integral average.
To this end, we first observe that $u$ vanishes continuously on
$\Delta_{2r}$.  Indeed, in the case that the continuous Dirichlet problem is solvable for $L$, this fact holds by assumption.
On the other hand,  in the case that $L$ is the heat operator, or an operator with $C^1$-Dini coefficients, 
we may use
the time-backwards ADR assumption (and the Wiener-type
criterion of \cite{EG}, or of \cite{FGL}) to deduce that
every point in $\Delta_{2r}$ is regular (see Remark \ref{r1.22}).
Thus, in either case, $u$ vanishes continuously on $\Delta_{2r}$.
We may then extend $u\equiv 0$ in $Q_{2r}\setminus \overline{\Omega}$, and we 
call this extension $\hat{u}$. Observe that $\hat{u}\geq 0$ and $\hat{u}$ is a subsolution in $Q_{2r}$. 
Therefore, by local boundedness \cite[Theorem 3]{M}, recalling that $T=T_{max}(Q_r)=t_0+r^2$,
we obtain
\begin{align} \label{eqB4}
\m(u) =
 \sup_{Q_{3r/2}\cap E(T)} \hat{u} 
\lesssim 
\tiltfiint_{Q_{2r}\cap E(T)}\hat{u} = Cr^{-n-2} \iint_{\Omega_{2r}(T)} u\,.
\end{align}
We note that Theorem 3 in \cite{M} is stated with an $L^p$ average, $p\geq 2$, on the 
right hand side of the inequality, but in hindsight,
this may be sharpened to an $L^1$ average, using a well-known self-improving property
of weak reverse H\"older estimates.  
\end{proof}

 \section{Proof of Lemma \ref{essbclosed}}\label{s6} 
 
 \begin{proof}[Proof of Lemma \ref{essbclosed}]
 We prove that the essential boundary $\eo$ and quasi-lateral boundary
 $\Sigma$ are closed sets.   In the case of $\eo$, by definition it is equivalent to show that
 the singular boundary $\so$ is relatively open in $\pom$.  To this end, fix
 ${\bf x_0}=(x_0,t_0)\in \so$, and note that by definition, there is an $\eps>0$ such that
 \[ Q_\eps^-:=Q_\eps^-({\bf x_0}) \subset \Omega\,, \quad \text{and}  \quad  Q_\eps^+:=Q_\eps^+({\bf x_0}) \subset \ree\setminus\Omega\,.\]
Since $Q_\eps^+$ is open, we therefore have that
\[Q_\eps^+ \subset \text{int} (\ree\setminus\Omega) =: \Omega_{ext}\,,\]
where int($A$) denotes the interior of $A$.  Consequently, if ${\bf x} = (x,t) \in \pom \cap Q_\eps$,
where $Q_\eps:= Q_\eps({\bf x_0})$, then
${\bf x}$ lies on the interface between $Q_\eps^-$ and $Q_\eps^+$, i.e., on the time-slice
$(Q_\eps)_{t_0}$. It follows that for any such ${\bf x}$, there is an $\eps' =\eps'({\bf x})>0$
such that
 \[ Q_{\eps'}^-({\bf x}) \subset \Omega\,, \quad \text{and}  \quad  Q_{\eps'}^+({\bf x}) \subset \ree\setminus\Omega\,,\]
i.e., ${\bf x}\in \so$, by definition, and thus $\so$ is relatively open, as desired.

To see that $\Sigma$ is closed, we first note that $(\so)_{T_{max}}$ is relatively open in $\po$, by the 
preceeding argument.
Thus, we need only observe in addition that the time slice
$(\bo)_{T_{min}}$ is also relatively open in $\pom$ 
(assuming that $T_{min}>-\infty$; otherwise there is nothing to prove).  
But this follows directly from the fact that
under the change of variable $t \to -t$, which maps $\Omega$ into an open set that we denote
$\Omega^*$, the time-slice $(\bo)_{T_{min}}$ is mapped onto $(\so^*)_{T_{max}(\Omega^*)}$. 
\end{proof}

\medskip
 
 \noindent {\bf Acknowledgements}.  
We thank Nicola Garofalo for bringing to our attention
 the references \cite{BiM} and \cite{GZ}, and Seick Kim for informing us of the results in \cite{BHL}, \cite{CDK},
 \cite{DK}, \cite{QX}, and \cite{SSSZ}.

\end{document}